\documentclass[11pt,reqno]{amsart}

\usepackage[text={160mm,240mm},centering]{geometry}            
\usepackage{amssymb,amsmath,setspace,geometry,indentfirst,changepage,inputenc,amsthm}
\usepackage{mathrsfs,amsfonts,savesym,graphicx,bm,color}

\usepackage{graphicx}
\usepackage{amssymb}
\usepackage{dsfont}

\usepackage{lscape}

\usepackage{tikz}
\usetikzlibrary{positioning}

\usepackage{float}

\usepackage[all,2cell]{xy}    
\UseAllTwocells               

\usepackage{amsmath,amsthm}

\usepackage[mathscr]{eucal}
\usepackage[all]{xy}
\usepackage{mathrsfs}
\usepackage{hyperref}
\usepackage{color}

\newtheorem{theorem}{Theorem}[section]
\newtheorem{lemma}[theorem]{Lemma}

\newtheorem{conjecture}[theorem]{Conjecture}
\newtheorem{corollary}[theorem]{Corollary}
\newtheorem{proposition}[theorem]{Proposition}

\newtheorem{definition-lemma}[theorem]{Definition-Lemma}
\newtheorem{definition-theorem}[theorem]{Definition-Theorem}

\theoremstyle{definition}
\newtheorem{example}[theorem]{Example}
\newtheorem{definition}[theorem]{Definition}

\newtheorem{remark}[theorem]{Remark}

\usepackage{fancyhdr}
\allowdisplaybreaks

\title{On the Monodromy Conjecture for Determinantal Varieties}

\author{Yifan Chen}
\address{
	Department of Mathematical Sciences,
	Tsinghua University,
	Beijing, 100084, P. R. China.}
\email{c-yf20@tsinghua.org.cn}

\author{Huaiqing Zuo}
\address{Department of Mathematical Sciences,
	Tsinghua University,
	Beijing, 100084, P. R. China.}
\email{hqzuo@mail.tsinghua.edu.cn}

\begin{document}
	
	\maketitle
	
	%
	
	\begin{abstract}
	This paper presents a proof of the monodromy conjecture for determinantal varieties. Our strategy centers on an in-depth analysis of monodromy zeta functions, leveraging a generalized A'Campo formula, an examination of multiple contact loci, and the exploitation of the intrinsic symmetric structures inherent to these varieties. Furthermore, we prove the holomorphy conjecture for determinantal varieties and the monodromy conjecture for Brill-Noether loci of generic curves.
	
	
	\end{abstract}
	
	\tableofcontents

	\section{Introduction}
	The monodromy conjecture stands as one of the most profound and enigmatic problems in modern singularity theory, residing at the confluence of number theory, algebra, analysis, geometry, and topology. It proposes a remarkable and deep connection between the arithmetic properties of a polynomial with integer coefficients and its geometric and topological invariants. 

	More explicitly, let $f \in \mathbb{Z}[x_{1},\dots,x_{n}]$ be a polynomial. On the arithmetic side, one associates to $f$ a $p$-adic Igusa zeta function $Z_{p,f}$, which encodes the number of solutions to $f \equiv 0 \mod p^{i}$ for all $i$. On the complex geometric side, by viewing $f$ as a complex polynomial, one studies the monodromy action on the cohomology of the local Milnor fibres of $f$. The conjecture boldly predicts that if $s_{0}$ is a pole of $Z_{p,f}$, then $e^{2\pi i \operatorname{Re}(s_{0})}$ must be an eigenvalue of this monodromy action.

	The framework was vastly generalized by Denef and Loeser, who introduced topological and motivic zeta functions, showing that the $p$-adic zeta function is a specialization of the universal motivic zeta function \cite{Zeta_functions, Germs_of_arcs_on_singular_varieties, Motivic_Igusa_Zeta_Function}. This allows one to formulate an analogous motivic monodromy conjecture. Furthermore, the conjecture has been extended to settings involving multiple polynomials or ideals \cite{Bernstein_Sato_polynomial_of_arbitrary_varieties, Proximite_evanescente_D_module}.

	Despite its fundamental importance, a general explanation for why the conjecture holds remains elusive. Known partial results are typically established on a case-by-case basis and require substantial, specific information about both the two aspects of the singularity. The reader is referred to \cite{Veys25} for an excellent survey on the monodromy conjecture.  In this paper, we provide a comprehensive proof of the monodromy conjecture for all determinantal varieties. Moreover, in a subsequent paper, we prove the monodromy conjecture for Pfaffians of a generic skew-symmetric matrix, a proof that relies partly on the method developed here.
	
	\vspace{0.5em}
	\noindent\textbf{Theorem A} (Theorem \ref{thmA})\textbf{.}
	\textit{The monodromy conjecture holds for determinantal varieties.}
	\vspace{0.5em}

	As an immediate corollary, the monodromy conjecture holds for the classical Brill-Noether loci of generic curves since locally they are determinantal varieties. Moreover, if the Petri map is injective, this result holds for higher rank vector bundles and their twists.  To be more precise, let $C$ be a smooth projective curve over $\mathbb{C}$ and let $F$ be a vector bundle on it. For integer $n \in \mathbb{Z}_{>0}$ and $d \in \mathbb{Z}_{\ge 0}$, we define $\mathcal{M}_{n,d}$ to be the moduli space of stable vector bundles on $C$ of rank $n$ and degree $d$. For $k \in \mathbb{Z}_{>0}$, let $\mathcal{V}_{n,d,k}(F):=\{E \in \mathcal{M}_{n,d}|h^{0}(C,E \otimes F) \ge k\}$ and it can be viewed as a closed subscheme of $\mathcal{M}_{n,d}$. Note that when $n=1$ and $F=\mathcal{O}_{C}$, $\mathcal{V}_{n,d,k}(F)$ is the classical Brill-Noether locus of $C$. For $E \in \mathcal{M}_{n,d}$, the Petri map is defined to be
	\begin{flalign*}
	\pi_{E,F}:H^{0}(C,E \otimes F) \otimes H^{0}(C, \check{E} \otimes \check{F} \otimes \omega_{C}) \rightarrow H^{0}(C,E \otimes \check{E} \otimes \omega_{C}),
	\end{flalign*}
	where $\check{E},\check{F}$ are the dual of $E,F$ and $\omega_{C}$ is the canonical sheaf of $C$. By Theorem $1.4$ in \cite{Budur_Brill-Noether_loci}, locally at $E$, $\mathcal{V}_{n,d,k}(F) \subset \mathcal{M}_{n,d}$ is defined by determinantal varieties. Combining with Theorem $1.1$, Theorem $1.5$ in \cite{Budur_Brill-Noether_loci} and Theorem A, we have the following corollary.
	\begin{corollary}\label{Budur_corollary}
	The monodromy conjecture holds locally at $E$ for the Brill-Noether loci $\mathcal{V}_{1,d,k}(\mathcal{O}_{C})$. In general, if the Petri map is injective, the monodromy conjecture holds locally at $E$ for $\mathcal{V}_{n,d,k}(F)$. 
	\end{corollary}

	Let $X=\mathbb{C}^{mn}(m \le n)$ be the space of $m \times n$ matrices and $Z_{r} \subset X$ be the subvariety of $X$ defined by $r \times r$ minors of matrices. $Z_{r}$ is called the determinantal variety. These varieties have many symmetric properties and they are good examples to compute for some important invariants of singularities. For example, in \cite{Arcs_on_Determinantal_Variteties_Roi}, Docampo analyzed the structure of the arc space of determinantal varieties, while in \cite{Stringy_E_function_of_determinantal_varieties}, the authors computed their stringy E-functions. The monodromy conjecture for determinantal varieties has been established for the case of maximal minors in \cite{Bernstein_for_maximal_minors} through the computation of the Bernstein–Sato polynomial of $Z_m$; however, our result addresses the general case and is therefore novel.
	
	One aspect of the monodromy conjecture—namely, the motivic zeta function of determinantal varieties—was computed by Docampo in \cite{Arcs_on_Determinantal_Variteties_Roi}. Furthermore, he explicitly determined the set of poles of the topological zeta function for $Z_r$, as in the following theorem. 
	
	\begin{theorem}[Theorem $6.5$ in \cite{Arcs_on_Determinantal_Variteties_Roi}]\label{roi_topological_zeta_function}
	With the notations above, when $m=n$, the topological zeta function of $Z_{r}$ is given by
	\begin{flalign*}
		Z^{\mathrm{top}}_{Z_{r}}(s)=\prod_{j=1}^{r}\frac{1}{1+s\cdot \frac{r+1-j}{(m+1-j)^{2}}}.
	\end{flalign*}
	\end{theorem}
	
	\begin{remark}\label{general_roi_zeta_function}
	Although this theorem is stated for the case when $m = n$, the same argument can be extended to general $m$ and $n$, yielding that $-\frac{(m+1-j)(n+1-j)}{r+1-j}(1 \le j \le r)$ are the poles of $Z^{\mathrm{top}}_{Z_{r}}(s)$.
	\end{remark}
	
	The present work focuses on the other side of the monodromy conjecture. In \cite{Denef_monodromy_zeta_function_and_eigenvalue} Denef proved that the zeros and poles of the monodromy zeta function can fully characterize the eigenvalues of monodromy, and Veys generalized this result to the case for ideals in \cite{Veys_generalized_monodromy_zeta_function}. This generalized monodromy zeta function, denoted by $Z^{\mathrm{mon}}_{Z_{r},e}(t)$, is defined over a point $e \in E$, where $E$ is the exceptional divisor of the blow-up of $Z_{r}$ in $X$. The proof of our main Theorem A is derived from a crucial result concerning the explicit computation of the monodromy zeta function for determinantal varieties, presented below.
	
	\vspace{0.5em}
	\noindent\textbf{Theorem B} (Theorem \ref{thmB})\textbf{.} 
	\textit{Let $\hat{X} \rightarrow X$ be the blow-up of $X$ at $Z_{r}$ and let $E$ be its exceptional divisor. For different points $e \in E$, the monodromy zeta function of $Z_{r}$ at $e$, denoted by $Z^{\mathrm{mon}}_{Z_{r},e}(t)$, can be $1$ or $1-t^{r+1-i}$ for $1 \le i \le r$, depending on the point $e$.}
	\vspace{0.5em}
	
	The core challenge in addressing this problem lies in characterizing the point $e \in E$. Because the blow-up of $Z_{r}$ in $X$ is a singular space with limited desirable properties, we instead analyze $e$ by transferring its information to the base space $X$. Specifically, we demonstrate that a non-trivial monodromy zeta function $Z^{\mathrm{mon}}_{Z_{r},e}(t)$ implies that $e$ can be characterized by a certain jet of $X$. We then classify the jets of $X$ that contain $e$ and compute their classes in the Grothendieck ring. These distinct classes correspond to the various possible forms of the monodromy zeta function, as presented in Theorem B.
	
	As an application of the tools and results developed in this paper, we prove the holomorphy conjecture for determinantal varieties. This conjecture, which is an analogue of the monodromy conjecture, gives a relation between the poles of twisted topological zeta function and the eigenvalues of the monodromy action, see Conjecture $3.4.2$ in \cite{Geometry_on_arc_spaces_of_algebraic_varieties}. In \cite{The_holomorphy_conjecture_for_ideals_in_dimension_two} this conjecture is generalized to the case of ideals. Details are provided in Section \ref{5}.
	
	\vspace{0.5em}
	\noindent\textbf{Theorem C} (Theorem \ref{thmC})\textbf{.}
	\textit{The holomorphy conjecture holds for determinantal varieties.}
	\vspace{0.5em}

	This paper is structured as follows. In Section \ref{sec2}, we provide a review of the generalized monodromy conjecture and foundational concepts such as jet schemes, arc spaces, and motivic zeta functions. Section \ref{sec3} is devoted to recalling key results concerning determinantal varieties. In Section \ref{sec4}, we give the proofs of Theorem A and Theorem B in three steps: we first connect the problem with multiple contact loci, then classify the required orbits and finally compute the monodromy zeta function. The holomorphy conjecture for determinantal varieties are introduced and proved in Section \ref{sec5}.
	
	\vspace{0.5em}
	\noindent\textbf{Acknowledgment}. We thank N. Budur for comments that led to Corollary \ref{Budur_corollary}. Y. Chen and H. Zuo were supported by BJNSF Grant 1252009. H. Zuo was supported by NSFC Grant 12271280.

	\section{Generalized monodromy conjecture}\label{sec2}
	\subsection{Jet schemes, the Grothendieck ring and motivic zeta functions}
	To introduce the monodromy conjecture more explicitly, we need the notion of motivic zeta function. We will review the basic concepts of jet schemes, arc spaces and the Grothendieck ring, using them to define motivic zeta functions. In this subsection we work over an algebraically closed field $k$ with characteristic $0$. 
	\begin{definition}[Jet Scheme]
		Let $X$ be a scheme over a field $k$, for any $m \in \mathbb{N}$, consider the functor from $k$-schemes to set
		\begin{flalign*}
			Z \mapsto \mathrm{Hom}(Z \times_{\mathrm{Spec}(k)} \mathrm{Spec} (k[t]/(t^{m+1})),X).
		\end{flalign*}
		There is a $k$-scheme representing this functor called the $m$-th jet scheme of $X$, and it is denoted by $J_{m}(X)$ (see Theorem 2.1, \cite{Jet_Schemes_Ishii}), i.e.,
		\begin{flalign*}
			\mathrm{Hom}(Z,J_{m}(X))=\mathrm{Hom}(Z \times_{\mathrm{Spec}(k)} \mathrm{Spec} (k[t]/(t^{m+1})),X).
		\end{flalign*}
	\end{definition}
	
	For $1 \le i \le j$, the truncation map $k[t]/(t^{j}) \to k[t]/(t^{i})$ induces natural projections $\psi_{i,j}: J_{j}(X) \to J_{i}(X)$ between jet schemes. If we identify $J_{0}(X)$ with $X$, we obtain a natural projection $\pi_{m}: J_{m}(X) \to X$. One can check that $\{J_{m}(X)\}_{m}$ is an inverse system and the inverse limit $J_{\infty}(X):= \varprojlim_m J_{m}(X)$ is called the arc space of $X$, inducing natural projections $\psi_{i}: J_{\infty}(X) \to J_{i}(X)$. If $X$ is of finite type, $J_{m}(X)$ is also of finite type for each $m \in \mathbb{N}$, but $J_{\infty}(X)$ is usually not.
	
	Next we introduce the Grothendieck ring, which makes the set of all algebraic varieties an abelian group.
	\begin{definition}[Grothendieck Ring]\label{Ch2_Sec2_Grothendieck_Ring}
		Let $k$ be a field and let $\mathrm{Var}_{k}$ be the category of $k$-varieties. The Grothendieck group of $k$-varieties, $K_0[\mathrm{Var}_k]$, is defined to be the quotient group of the free abelian group with basis $\{[X]\}_{X\in\mathrm{Var}_k}$, modulo the following relations.
		\begin{align*}
			& [X]-[Y], \ \mathrm{if} \ X \simeq Y,\\
			& [X]-[X_{\mathrm{red}}],\\
			& [X]-[U]-[X\setminus U],\ \mathrm{for} \ \mathrm{any} \ \mathrm{open} \ \mathrm{set} \ U\subseteq X.\ 
		\end{align*}
		One can further define a multiplication structure on $K_0[\mathrm{Var}_k]$ by setting
		\begin{displaymath}
			[X] \cdot [Y] := [X\times Y].
		\end{displaymath}
		This makes $K_0[\mathrm{Var}_k]$ a ring, called the Grothendieck ring of $k$-varieties. In particular, let $\mathbb L := [\mathbb A_k^1]$ and $K_0[\mathrm{Var}_k]_{\mathbb L}$ be the localization at $\mathbb L$.
	\end{definition}
	
	\begin{remark}
	It follows from the definition that the Euler characteristic $\chi: \mathrm{Var}_{k} \rightarrow \mathbb{Z}$ factors through the Grothendieck ring, therefore we obtain the Euler characteristic specialization $\chi: K_{0}[\mathrm{Var}_{k}] \rightarrow \mathbb{Z}$.
	\end{remark}

	\noindent\textbf{Kontsevich's Completion}
	
	Kontsevich's completion allows us to take limits in Grothendieck ring. To introduce this, we define a decreasing filtration $F^\bullet$ on $K_0[\mathrm{Var}_{k}]_{\mathbb L}$ as follows. For $m\in \mathbb Z$, $F^m$ is the subgroup of $K_0[\mathrm{Var}_{k}]_{\mathbb L}$ generated by $[S] \cdot \mathbb L^{-i}$ with $\dim S- i \leq -m$. It is actually a ring filtration, that is, $F^m \cdot F^n \subseteq F^{m+n}$.
	\begin{definition}
		The Kontsevich's completed Grothendieck ring $\widehat{K_0[\mathrm{Var}_{k}]}$ is defined to be 
		\begin{displaymath}
			\widehat{K_0[\mathrm{Var}_{k}]} := \varprojlim_{m\in \mathbb Z} K_0[\mathrm{Var}_{k}]_{\mathbb L}/F^m.
		\end{displaymath}
	\end{definition}

	Now we fix an ideal $I \subset k[x_{1},\dots,x_{n}]$, and we will give the definition of motivic zeta function for $I$. Before that, we need the definition of contact loci.
	
	\begin{definition}[Contact locus]
	For any polynomial $f \in k[x_{1},\dots,x_{n}]$ and any arc $\psi \in J_{\infty}(\mathbb{A}_{k}^{n})$ or any jet $\psi \in J_{m}(\mathbb{A}_{k}^{n})$ for some $m$, the order of $f$ along $\psi$, denoted by $\mathrm{ord}_{f}(\psi)$, is defined to be the order of the formal power series $f(\psi)$. For any ideal $J \subset k[x_{1},\dots,x_{n}]$, the order of $J$ along $\psi$, denoted by $\mathrm{ord}_{J}(\psi)$, is the minimum of $\mathrm{ord}_{f}(\psi)$ for all $f \in J$. The $m$-th contact locus of $J$ is defined to be
	\begin{flalign*}
		\mathfrak{X}_{m,J}:=\{\psi \in J_{\infty}(\mathbb{A}_{k}^{n})\; | \;\mathrm{ord}_{J}(\psi)=m\}.
	\end{flalign*}
	Its truncation to the $m$-jet level is 
	\begin{flalign*}
		\mathfrak{X}_{m,J}^{m}:=\{\psi \in J_{m}(\mathbb{A}_{k}^{n})\; | \;\mathrm{ord}_{J}(\psi)=m\}.
	\end{flalign*}
	\end{definition}
	
	Now we can define the motivic zeta function for the ideal $I$.
	\begin{definition}[Motivic zeta function]
	For an ideal $I \subset k[x_{1},\dots,x_{n}]$, the motivic zeta function of $I$ is given by
	\begin{flalign*}
		Z_{I}^{\mathrm{mot}}(s):=\sum_{m=0}^{\infty}[\mathfrak{X}_{m,I}^{m}]\mathbb{L}^{-ms-mn}.
	\end{flalign*}
	In particular, if $I$ is generated by one polynomial $f$, then the motivic zeta function of $f$ is defined to be
	\begin{flalign*}
		Z_{f}^{\mathrm{mot}}(s):=\sum_{m=0}^{\infty}[\mathfrak{X}_{m,(f)}^{m}]\mathbb{L}^{-ms-mn}.
	\end{flalign*}
	\end{definition}
	
	By definition, the motivic zeta function lies in $K_0[\mathrm{Var}_{k}]_{\mathbb L}$. Denef and Loeser showed that the motivic zeta function can be computed by log resolution and it is a "rational" function. Although their version is about one polynomial, the result also holds for ideals.
	\begin{theorem}[\cite{Motivic_Igusa_Zeta_Function}]\label{theorem_denef_loeser}
	Let $\mu:Y \rightarrow \mathbb{A}_{k}^{n}$ be a log resolution of $(\mathbb{A}_{k}^{n},V(I))$. Suppose $E_{i}(i \in S)$ are the irreducible components of $\mu^{*}(I)$ and relative canonical divisor $K_{Y/\mathbb{A}_{k}^{n}}$ such that $\mu^{*}(I)=\sum_{i \in S}N_{i}E_{i}$ and $K_{Y/\mathbb{A}_{k}^{n}}=\sum_{i \in S}(\nu_{i}-1)E_{i}$, then
	\begin{flalign*}
		Z^{\mathrm{mot}}_{I}(s)=\sum_{J \subset S}[\mathring{E}_{J}]\prod_{i \in I}\frac{\mathbb{L}-1}{\mathbb{L}^{N_{i}s+\nu_{i}}-1},
	\end{flalign*}
	where $\mathring{E}_{J}:=(\cap_{i \in J}E_{i})\setminus (\cup_{j \notin J}E_{j})$. Here we call $(N_{i},\nu_{i})(i \in S)$ the data of the log resolution.
	\end{theorem} 
	
	The motivic zeta function can be specialized into the topological zeta function by imposing the Euler characteristic $\chi(\cdot)$ into the expression, as in the following definition.
	\begin{definition}[Topological zeta function]
	As in the setting above, the topological zeta function of $I$ is given by
	\begin{flalign*}
		Z^{\mathrm{top}}_{I}(s):=\sum_{J \subset S}\chi(\mathring{E}_{J})\prod_{i \in J}\frac{1}{N_{i}s+\nu_{i}}.
	\end{flalign*}
	\end{definition}
	
	From this expression, it follows that the numbers $-\frac{\nu_{i}}{N_{i}}$ are the "candidate poles" of the motivic and topological zeta function. In practice, people often find that some of the candidate poles are cancelled. Understanding these poles is one of the major problems of the monodromy conjecture.  
	
	\subsection{Monodromy conjecture}\label{monodromy_conj}
	
	The monodromy conjecture, which was originally raised for one polynomial, establishes the relation between the poles of the motivic zeta function, the eigenvalues of local Milnor monodromy action and the roots of Bernstein-Sato polynomials. Given non-constant polynomial $f \in \mathbb{C}[x_{1},\dots,x_{n}]$ and a point $x \in V(f)$, there is a natural monodromy action on the cohomology groups of the Milnor fibre of $f$ at $x$. We denote by $E(f)$ the set of all eigenvalues of the cohomology groups of the Milnor fiber at any point $x \in V(f)$.
	
	\begin{conjecture}[Monodromy Conjecture, \cite{Motivic_Igusa_Zeta_Function}]\label{mc}
	Let $f\in \mathbb C[x_{1},\dots,x_{n}]\setminus \mathbb C$. If $s_0$ is a pole of $Z_{f}^{\mathrm{mot}}(s)$, then $e^{s_0\cdot 2\pi \sqrt{-1}} \in E(f)$.
	\end{conjecture}

	\begin{remark}
	There is another stronger version of monodromy conjecture which states that the poles of motivic zeta function are the roots of Bernstein-Sato polynomial of $f$. In fact, the authors of \cite{Kashiwara_Bernstein_and_eigenvalue_of_monodromy} and \cite{Malgrange_Bernstein_eigenvalue} showed that the set $$\mathrm{Exp}(V(b_{f})):=\{e^{2\pi \sqrt{-1}s_{0}}\; | \; b_{f}(s_{0})=0\}$$ equals $E(f)$, therefore this version is stronger than the original one. For the introduction of Bernstein-Sato polynomial, one can refer to, for example, \cite{Rationality_of_Roots_of_B-Function_Kashiwara}, \cite{Budur_Bernstein_notes}, and \cite{Popa_notes}. 
	\end{remark}

	Next we introduce another kind of zeta function, called the monodromy zeta function $Z^{\mathrm{mon}}_{f, e}$, which is defined over $f$ and a point $e$ of $V(f)$ and has a close relationship with the eigenvalue of the monodromy action. The definition is as follows. 
	
	\begin{definition}
	For $f \in \mathbb{C}[x_{1},\dots,x_{n}] \setminus \mathbb{C}$ and $x \in V(f)$, the monodromy zeta function of $f$ at $x$ is defined as
	\begin{flalign*}
		Z^{\mathrm{mon}}_{f, x}:=\prod_{j=0}^{n-1}P_{j}(t)^{(-1)^{j-1}},
	\end{flalign*}
	where $P_{j}(t)$ is the characteristic polynomial of the monodromy action on $H^{i}(F_{f, x},\mathbb{C})$ and $F_{f, x}$ is the Milnor fibre of $f$ at $x$.
	\end{definition}

	The monodromy zeta function can be computed with log resolution using the following A'Campo formula.

	\begin{theorem}[A'Campo formula]
	Let $\mu:Y \rightarrow \mathbb{C}^{n}$ be a log resolution of $(\mathbb{C}^{n},V(f))$. Suppose $E_{i}(i \in S)$ are the irreducible components of $\mu^{*}(V(f))$ such that $\mu^{*}(V(f))=\sum_{i \in S}N_{i}E_{i}$, then for any point $e \in V(f)$, the monodromy zeta function of $f$ at $e$ is given by
	\begin{flalign*}
		Z^{\mathrm{mon}}_{f, e}(t)=\prod_{i \in S}(1-t^{N_{i}})^{\chi(\mathring{E_{i}} \cap \mu^{-1}(e))},
	\end{flalign*}
	where $\chi( \cdot )$ is the Euler characteristic.
	\end{theorem}
	
	\begin{remark}\label{little_remark}
	It is shown in \cite{Denef_monodromy_zeta_function_and_eigenvalue} that the union of zeros and poles of the monodromy zeta function of all points in $V(f)$ equals $E(f)$, therefore the monodromy zeta function can fully characterize $E(f)$. 
	\end{remark}
	
	The monodromy conjecture can be generalized to the ideal case. When we replace $f$ with an ideal $I \subset \mathbb{C}[x_{1},\dots,x_{n}]$, we can replace $Z^{\mathrm{mot}}_{f}$ with $Z^{\mathrm{mot}}_{I}$. By \cite{Bernstein_Sato_polynomial_of_arbitrary_varieties}, the Bernstein-Sato polynomial can also be generalized to the ideal case, denoted by $b_{I}(s)$. Unfortunately there is no counterpart of the local Milnor monodromy eigenvalue for ideal version. However, Verdier introduced the notion of Verdier monodromy in \cite{Verdier_new_monodromy} which is equivalent to the Milnor monodromy version and it can be naturally generalized to the ideal case. This generalization is compatible with that of the Bernstein-Sato polynomial, since in \cite{Bernstein_Sato_polynomial_of_arbitrary_varieties} Budur, Musta\c{t}ă and Saito proved that the set $$\mathrm{Exp}(V(b_{I}(s))):=\{e^{2\pi \sqrt{-1}s_{0}}\; | \; b_{I}(s_{0})=0\}$$ equals the set of the union of all eigenvalues of Verdier monodromy. For the monodromy zeta function part, in \cite{Veys_generalized_monodromy_zeta_function}, Proeyen and Veys defined monodromy zeta function for ideals using Verdier monodromy. This generalization is also compatible with other parts as the following proposition.
	\begin{proposition}[Remark $3.1$ in \cite{Veys_generalized_monodromy_zeta_function}]\label{ideal_case_monodromy_zeta_function_and_eigenvalue}
	The union of all zeros and poles of monodromy zeta function for all possible points is equal to the set of all eigenvalues of Verdier monodromy.
	\end{proposition}

	 They also generalized A'Campo formula as follows.
	
	\begin{theorem}[Theorem $3.2$ in \cite{Veys_generalized_monodromy_zeta_function}]\label{generalized_monodromy_zeta_func_veys}
	Let $I \subset \mathbb{C}[x_{1},\dots,x_{n}]$ be an ideal and $h: \hat{\mathbb{C}}^{n} \rightarrow \mathbb{C}^{n}$ is the blow-up of $I$. Suppose $\mu: Y \rightarrow \mathbb{C}^{n}$ be the log resolution of $(\mathbb{C}^{n},V(I))$ and $\phi: Y \rightarrow \hat{\mathbb{C}}^{n}$ be the morphism such that $h \circ \phi = \mu$. Assume $E_{i}(i \in S)$ are the irreducible components of $\mu^{*}(I)$ and $\mu^{*}(I)=\sum_{i \in S}N_{i}E_{i}$, then for any point $e \in h^{-1}(V(I))$, the monodromy zeta function of $I$ at $e$ is given by
	\begin{flalign*}
		Z^{\mathrm{mon}}_{V(I),e}=\prod_{i \in S}(1-t^{N_{i}})^{-\chi(\mathring{E_{i}} \cap \phi^{-1}(e))}.
	\end{flalign*}
	\end{theorem}

	\section{Results about determinantal varieties}\label{sec3}
	
	Before starting the proofs for Theorem A and Theorem B, we establish notations and review fundamental results concerning determinantal varieties and their arc spaces. We also introduce the result about determinantal variety defined by maximal minors in \cite{Bernstein_for_maximal_minors}.
	
	Let $X=\mathbb{C}^{mn}$$(m \le n)$ be the space of $m \times n$ matrices and $Z_{r}$$(1 \le r \le m)$ be the scheme of matrices with rank $\le r-1$, so $Z_{r}$ is defined by the ideal of $r$-minors, which we denote by $I_{r}$. Let $x_{i,j}(1 \le i \le m, 1 \le j \le n)$ be the elements of the matrix, and the coordinate ring of $X$ is $\mathbb{C}[x_{i,j}\;|\; 1 \le i \le m, 1 \le j \le n]$. Let $G=\mathrm{GL}_{m} \times \mathrm{GL}_{n}$ and we define the action of $G$ on X as follows.
	\begin{flalign*}
		\mathrm{GL}_{m} \times \mathrm{GL}_{n} \times X &\longrightarrow X   \\
		(g,h) \cdot A &\mapsto gAh^{-1}.
	\end{flalign*}
	This action induces the action of $J_{\infty}(G)$ on $J_{\infty}(X)$ by just replacing the elements in matrices by the formal power series. Now $J_{\infty}(X)$ is decomposed into some orbits by this action, and every orbit has a standard form, as the following theorem in \cite{Arcs_on_Determinantal_Variteties_Roi}.
	\begin{theorem}[Proposition 3.2 in \cite{Arcs_on_Determinantal_Variteties_Roi}]\label{Orbit decomposition}
		Suppose $0 \le \lambda_{1} \le \lambda_{2} \le \dots \le \lambda_{m} \le \infty$ with $\lambda_{1},\dots,\lambda_{m} \in \mathbb{N} \cup \{\infty\}$, and $\lambda=(\lambda_{1},\dots,\lambda_{m})$. Then every standard form of the orbit is of the form $\delta_{\lambda}$, defined as follows:
		\begin{flalign*}
			\delta_{\lambda}=\begin{pmatrix}
				t^{\lambda_{1}} & 0 & \dots & 0 & \dots & 0 \\
				0 & t^{\lambda_{2}} & \dots & 0 & \dots & 0 \\
				\vdots & \vdots & \ddots & \vdots & \dots & 0 \\
				0 & 0 & \dots & t^{\lambda_{m}} & \dots & 0
			\end{pmatrix}.
		\end{flalign*}
		Here we set $t^{\infty}=0$. Thus every orbit is indexed by some $\lambda$ and we denote the orbit corresponding to $\lambda$ by $\mathcal{C}_{\lambda}$. Moreover, $J_{\infty}(Z_{r})$ is stable under the $J_{\infty}(G)$-action, and we have $\mathcal{C}_{\lambda} \subset J_{\infty}(Z_{r})$ if and only if $\lambda_{r}=\dots=\lambda_{m}=\infty$.
	\end{theorem}
	
	Similarly for any positive integer $l$ we have the action of $J_{l}(G)$ on $J_{l}(X)$. The orbit decomposition and the standard form of each orbit are similar to that of arc space, as given in the following proposition in \cite{Arcs_on_Determinantal_Variteties_Roi}.
	\begin{proposition}[Proposition $3.9$ in \cite{Arcs_on_Determinantal_Variteties_Roi}]
	For any $\lambda=(\lambda_{1},...,\lambda_{m})$ with $0 \le \lambda_{1} \le \dots \le \lambda_{m} \le \infty$, we define $\bar{\lambda}:=(\bar{\lambda}_{1},...,\bar{\lambda}_{m})$, where $\bar{\lambda}_{i}=\mathrm{min}\{l+1,\lambda_{i}\}$, then every standard form of the orbit is of the form $\delta_{\bar{\lambda}}$. We denote the orbit on the $l$-th jet level corresponding to $\delta_{\bar{\lambda}}$ by $\mathcal{C}_{\bar{\lambda}}^{l}$, then the inverse image of $\mathcal{C}_{\mu}^{l}$ under the truncation map $J_{\infty}(X) \rightarrow J_{l}(X)$ is the union of orbits $\mathcal{C}_{\lambda}$ such that $\bar{\lambda}=\mu$. 	
	\end{proposition}
	
	The log resolution of $(X,Z_{r})$ and the corresponding data $(N_{i},\nu_{i})$ are given in \cite{Multiplier_Ideals_of_Determinantal_Ideal}, with the following proposition.
	\begin{proposition}[Theorem $4.4$, Corollary $4.5$ and Corollary $4.6$ in \cite{Multiplier_Ideals_of_Determinantal_Ideal}]\label{log_resolution_of_determinantal_variety}
		Let $\pi_{1}: A_{1} \longrightarrow A_{0}:=\mathbb{C}^{mn}$ be the blowup of $Z_{0}$. Suppose $\pi_{1}:A_{1} \longrightarrow A_{0},\dots,\pi_{i-1}:A_{i-1} \longrightarrow A_{i-2}$ have been defined, set $\pi_{i}: A_{i} \longrightarrow A_{i-1}$ to be the blowup of $\tilde{Z}_{i-1}$ in $A_{i-1}$, where $\tilde{Z}_{i-1}$ is the strict transform of $Z_{i-1}$ via $\pi_{1} \circ \pi_{2} \circ \dots \circ \pi_{i-1}:A_{i-1} \longrightarrow A_{0}$. Then $\pi: A_{r} \longrightarrow A_{r-1} \longrightarrow \dots \longrightarrow A_{1} \longrightarrow \mathbb{C}^{mn}$ is a log resolution of $(X,Z_{r})$. Moreover, if we denote the irreducible component of the exceptional divisor of $\pi$ obtained from the exceptional divisor of $\pi_{j}$ by $E_{j}$$(1 \le j \le r)$, then we have $\pi^{*}(Z_{r})=\sum_{j=1}^{r}(r+1-j)E_{j}$ and $K_{A_{r}/A_{0}}=\sum_{j=1}^{r}((m+1-j)(n+1-j)-1)E_{j}$, i.e. the data $(N_{j},\nu_{j})=(r+1-j,(m+1-j)(n+1-j))$ for $1 \le j \le r$.
	\end{proposition}
	
	\begin{remark}\label{poles_of_zeta_function}
		From the proposition above we know $-\frac{(m+1-j)(n+1-j)}{r+1-j}$ $(1 \le j \le r)$ are the candidate poles of $Z^{\mathrm{mot}}_{Z_{r}}$.  The general case of Theorem \ref{roi_topological_zeta_function} (see Remark \ref{general_roi_zeta_function}) shows that these candidate poles are all actual poles of motivic zeta function. 
	\end{remark}
	
	When $m=r$, $Z_{r}$ is defined by the maximal minors of the matrix. In this case, the authors of \cite{Bernstein_for_maximal_minors} give the Bernstein-Sato polynomial of $I_{r}$, proving the monodromy conjecture for $I_{r}$. The result is as follows.
	\begin{theorem}[\cite{Bernstein_for_maximal_minors}]
		The Bernstein-Sato polynomial of $I_{r}$ is given by
		\begin{flalign*}
			b_{I_{r}}(s)=\prod_{j=n-m+1}^{n}(s+j).
		\end{flalign*}
		In particular, monodromy conjecture for $I_{r}$ holds since the candidate poles of $Z^{\mathrm{mot}}_{Z_{r}}(s)$ are $-(n+1-j)$ $(1 \le j \le m)$.
	\end{theorem}

	\section{Proofs of Theorem A and Theorem B}\label{sec4}
	
	We adopt the notation from the previous section. In order to prove the monodromy conjecture for determinantal varieties, we will calculate monodromy zeta function for the ideal $I_{r}$, whose zeros and poles will give the eigenvalues of the monodromy action. We will use Theorem \ref{generalized_monodromy_zeta_func_veys}, and firstly let us fix some notations. Suppose $h:\hat{X} \rightarrow X$ is the blow-up of $Z_{r}$ in $X$. Let $f_{1},\dots,f_{w}$ be the all $r \times r$ minors, and $I_{r}=(f_{1},\dots,f_{w})$, then $D_{+}(f_{i}) \subset \hat{X}=\mathrm{Proj}\bigoplus_{j \ge 0}I_{r}^{j}$ is an open subset for all $i$. Moreover, we have $D_{+}(f_{i}) \cong \mathrm{Spec} \mathbb{C}[x_{1,1},\dots,x_{m,n},\frac{f_{1}}{f_{i}},\dots,1,\dots,\frac{f_{w}}{f_{i}}]$. We denote $h^{-1}(Z_{r})$ by $E$ and it is defined by $f_{i}$ in $D_{+}(f_{i})$ for all $i$. For a point in $\hat{X}$, we can write $(a_{1,1},\dots,a_{m,n},[b_{1},\dots,b_{w}])$ as the coordinates, where $[ \cdot]$ denotes points in projective space, i.e. $[b_{1},\dots,b_{w}]=[\theta b_{1},\dots,\theta b_{w}]$ for any non-zero $\theta$. Let $\mu: Y \rightarrow X$ be the canonical resolution of $(X,Z_{r})$ introduced in Proposition \ref{log_resolution_of_determinantal_variety}, then by the universal property of blow-up, there exists a morphism $\phi: Y \rightarrow \hat{X}$ such that $\mu=h \circ \phi$.  Assume $E_{i}(i \in S)$ are the irreducible components of $\mu^{-1}(Z_{r})$, $\mu^{*}(Z_{r})=\sum_{i \in S}N_{i}E_{i}$ and $K_{Y/X}=\sum_{i \in S}(\nu_{i}-1)E_{i}$, then by Theorem \ref{generalized_monodromy_zeta_func_veys}, for $e \in E$, the monodromy zeta function of $I_{r}$ at $e$ is
	\begin{flalign*}
		Z^{\mathrm{mon}}_{Z_{r},e}(t)=\prod_{i \in S}(1-t^{N_{i}})^{\chi(\mathring{E_{i}} \cap \phi^{-1}(e))}.
	\end{flalign*}
	Actually by Proposition \ref{log_resolution_of_determinantal_variety}, we have $S=\{1,\dots,r\}$, $E_{i}$ is the strict transform of $Z_{i}$ in $Y$, $N_{i}=r+1-i$ and $\nu_{i}=(m+1-i)(n+1-i)$ for $i=1,\dots,r$. This implies that it suffices to calculate $\chi(\mathring{E_{i}} \cap \phi^{-1}(e))$ for every $1 \le i \le r$ and $e \in E$. 
	
	\vspace{0.5em}
	\noindent\textbf{Step 1: Connection with multiple contact loci}.
	\vspace{0.5em}
	
	We will connect this space with the contact loci of $X$ with respect to $E_{i}$. First we state the proposition which gives the relation of the arc spaces of $X$ and $Y$.
	\begin{proposition}[\cite{Jet_schemes_and_singularities}, Proposition $3.2$]
	Suppose $f:X' \rightarrow X$ is a proper morphism and $Z \subset X$ is a closed subset such that $f$ induces an isomorphism between $X'-f^{-1}(Z)$ and $X-Z$, then the morphism on the arc space level $f_{\infty}:J_{\infty}(X') \rightarrow J_{\infty}(X)$ induces a bijection of $J_{\infty}(X')-J_{\infty}(f^{-1}(Z))$ and $J_{\infty}(X)-J_{\infty}(Z)$.
	\end{proposition} 
	
	Now we take $Z=Z_{r}$ and $X'=Y$, then for every arc $\gamma \in J_{\infty}(X)-J_{\infty}(Z_{r})$, we can lift it into an arc $\tilde{\gamma} \in J_{\infty}(Y)$. For a tuple of integers $(u_{1},\dots,u_{r}) \in \mathbb{N}^{r}$, we define the multiple contact loci with respect to $E_{i}$ to be
	\begin{flalign*}
		\mathfrak{X}_{u_{1},\dots,u_{r}}:=\{\gamma \in J_{\infty}(X)-J_{\infty}(Z_{r})\; | \; \mathrm{ord}_{E_{i}}(\tilde{\gamma})=u_{i},1 \le i \le r\}.
	\end{flalign*}
	Its corresponding set in $J_{\infty}(Y)$ is exactly
	\begin{flalign*}
		\mathcal{Y}_{u_{1},\dots,u_{r}}:=\{\gamma \in J_{\infty}(Y)\; | \;\mathrm{ord}_{E_{i}}(\gamma)=u_{i},1 \le i \le r\}.
	\end{flalign*}
	Note that for $\gamma \in \mathcal{Y}_{u_{1},\dots,u_{r}}$, we have $\mathrm{ord}_{\mu^{-1}(Z_{r})}(\gamma)=N$, where $N=\sum_{i=1}^{r}N_{i}u_{i}$. Let $l$ be an integer big enough such that $l \gg N$, and we truncate the contact loci to the $l$-th jet level to obtain
	\begin{flalign*}
		\mathcal{Y}_{u_{1},\dots,u_{r}}^{l}:=\{\gamma \in J_{l}(Y)\; | \; \mathrm{ord}_{E_{i}}(\gamma)=u_{i},1 \le i \le r\}=\psi_{l}^{Y}(\mathcal{Y}_{u_{1},\dots,u_{r}})
	\end{flalign*}  
	and
	\begin{flalign*}
		\mathfrak{X}_{u_{1},\dots,u_{r}}^{l}:=\psi_{l}^{X}(\mathfrak{X}_{u_{1},\dots,u_{r}}),
	\end{flalign*}
	where $\psi_{l}^{X}:J_{\infty}(X) \rightarrow J_{l}(X)$ and $\psi_{l}^{Y}:J_{\infty}(Y) \rightarrow J_{l}(Y)$ are the natural truncation map. Suppose $I_{u_{1},\dots,u_{r}}=\{1 \le i \le r \; | \; u_{i} \neq 0\}$, then $\pi_{l}^{Y}(\mathcal{Y}_{u_{1},\dots,u_{r}}) \in \mathring{E}_{I_{u_{1},\dots,u_{r}}}$, where $\pi_{l}^{Y}:J_{l}(Y) \rightarrow Y$ is the projection map. With these notations, we introduce the following two propositions.
	\begin{proposition}[Theorem $3.4$ in \cite{cohomology_of_contact_loci}]\label{lemmadefloeser}
 	For $l$ large enough, $\mu_{l}:\mathcal{Y}^{l}_{u_{1},\dots,u_{r}} \rightarrow \mathfrak{X}^{l}_{u_{1},\dots,u_{r}}$ is a Zariski locally trivial fibration with fibre $\mathbb{C}^{\sum_{i=1}^{r}u_{i}(\nu_{i}-1)}$. 
	\end{proposition}
	\begin{proposition}[Proposition $3.2$ in \cite{on_the_embedded_nash_problem}]\label{lemmabudur}
	$\pi_{l}^{Y}|_{\mathcal{Y}^{l}_{u_{1},\dots,u_{r}}}:\mathcal{Y}^{l}_{u_{1},\dots,u_{r}} \rightarrow \mathring{E}_{I_{u_{1},\dots,u_{r}}}$ is a Zariski locally trivial fibration with fibre $(\mathbb{C}^{*})^{|I_{u_{1},\dots,u_{r}}|
	} \times \mathbb{C}^{mnl-\sum_{i=1}^{r}u_{i}}$. 
	\end{proposition}
	
	\begin{remark}
	Note that Proposition \ref{lemmabudur} was stated for the case when $|I_{u_{1},\dots,u_{r}}|=1$, but the result extends to the general case as above with the same proof. 
	\end{remark}

	Proposition \ref{lemmabudur} and Proposition \ref{lemmadefloeser} relate the Grothendieck class of $\mathring{E}_{I}$ to that of the multiple contact loci. For our use, we only need to take $u_{1}=\dots=u_{i-1}=u_{i+1}=\dots=u_{r}=0$ and $u_{i}=1$, then by the proposition above, $\pi_{l}^{Y}: (\pi_{l}^{Y})^{-1}(\mathring{E}_{i} \cap \phi^{-1}(e))\cap \mathcal{Y}^{l}_{0,\dots,1,\dots,0} \rightarrow \mathring{E}_{i} \cap \phi^{-1}(e)$ is a locally trivial fibration, therefore
	\begin{align}\label{1}
	[\mathring{E}_{i} \cap \phi^{-1}(e)]=\frac{[(\pi_{l}^{Y})^{-1}(\mathring{E}_{i} \cap \phi^{-1}(e))\cap \mathcal{Y}^{l}_{0,\dots,1,\dots,0}]}{(\mathbb{L}-1)\mathbb{L}^{mnl-1}}.
	\end{align}
	
	Next we will prove that $(\pi_{l}^{Y})^{-1}(\mathring{E}_{i} \cap \phi^{-1}(e))\cap \mathcal{Y}^{l}_{0,\dots,1,\dots,0}$ is also a union of fibres of the map $\mu_{l}$.
	\begin{proposition}\label{alsofibre}
	If $\mathring{E}_{i} \cap \phi^{-1}(e)$ is not empty, then $\mu_{l}:(\pi_{l}^{Y})^{-1}(\mathring{E}_{i} \cap \phi^{-1}(e))\cap \mathcal{Y}^{l}_{0,\dots,1,\dots,0} \rightarrow \mu_{l}((\pi_{l}^{Y})^{-1}(\mathring{E}_{i} \cap \phi^{-1}(e))\cap \mathcal{Y}^{l}_{0,\dots,1,\dots,0})$ is also a Zariski locally trivial fibration with fibre $\mathbb{C}^{\nu_{i}-1}$.
	\end{proposition}
	
	\begin{proof}
	We can rewrite the set on the left hand side as
	\begin{flalign*}
		&(\pi_{l}^{Y})^{-1}(\mathring{E}_{i} \cap \phi^{-1}(e))\cap \mathcal{Y}^{l}_{0,\dots,1,\dots,0}  \\
		=&\{\gamma \in J_{l}(Y) \; | \; \mathrm{ord}_{E_{i}}(\gamma)=1,\gamma(0) \in \mathring{E}_{i},\phi(\gamma(0))=e\}  \\
		=&\{\gamma \in J_{l}(Y)\; | \; \mathrm{ord}_{E_{i}}(\gamma)=1,\gamma(0) \in \mathring{E}_{i},\phi_{l}(\gamma)(0)=e\}.	
	\end{flalign*}
	Note that since the strict transform of $Z_{j}$ in $Y$ is $jE_{1}+(j-1)E_{2}+\dots+E_{j}$ for all $1 \le j \le r$, the condition of the order of $\gamma$ along each $E_{j}$ is equivalent to $\mu_{l}(\gamma) \in \mathcal{C}_{\lambda}^{l}$, where $\lambda=(\lambda_{1},\dots,\lambda_{m})$ satisfies that $\lambda_{1}=\dots=\lambda_{i-1}=0$ and $\lambda_{i}=\dots=\lambda_{r}=1$. We may assume $e \in D_{+}(f_{1})=\mathrm{Spec}\mathbb{C}[x_{1,1},\dots,x_{m,n},\frac{f_{2}}{f_{1}},\dots,\frac{f_{w}}{f_{1}}]$ and $e=(a_{1,1},\dots,a_{m,n},b_{2},\dots,b_{w})$ written in coordinates. We set $\hat{\gamma}:=h_{l}(\phi_{l}(\gamma))=\mu_{l}(\gamma) \in \mathcal{C}_{\lambda}^{l}$. Since $\hat{\gamma} \in \mathcal{C}_{\lambda}^{l}$, the condition $\phi_{l}(\gamma)(0)=e$ is equivalent to $\hat{\gamma}(0)=(a_{1,1},\dots,a_{m,n})$, $\mathrm{ord}_{t}(f_{1}(\hat{\gamma}))=\mathrm{min}\{f_{j}(\hat{\gamma})\}$ and $\frac{f_{j}(\hat{\gamma})}{f_{1}(\hat{\gamma})}(0)=b_{j}$ for $2 \le j \le w$. We denote the set in $J_{l}(X)$ satisfying these three conditions by $\mathcal{D}_{l,e}$, then we have
	\begin{align}\label{2}
		(\pi_{l}^{Y})^{-1}(\mathring{E}_{i} \cap \phi^{-1}(e))\cap \mathcal{Y}^{l}_{0,\dots,1,\dots,0}=\bigcup_{\lambda_{1}=\dots=\lambda_{i-1}=0,\lambda_{i}=\dots=\lambda_{r}=1}\mu_{l}^{-1}(\mathcal{C}_{\lambda}^{l} \cap \mathcal{D}_{l,e}).
	\end{align}
	Now the assertion follows from Proposition \ref{lemmadefloeser}.
	\end{proof}
	
	\begin{remark}\label{only_one_e_nonzero}
	By the proof of Proposition \ref{alsofibre}, when $\mathring{E_{i}} \cap \phi^{-1}(e)$ is not empty, there exists $\hat{\gamma} \in \mathcal{C}_{\lambda}^{l}$ such that $\hat{\gamma}(0)=h(e)$. This implies that the rank of $h(e)$ is $i-1$, therefore for any fixed $e \in E$, there is at most one $i$ such that $\chi(\mathring{E_{i}} \cap \phi^{-1}(e)) \neq 0$, depending on the rank of $h(e)$.
	\end{remark}

	\vspace{0.5em}
	\noindent\textbf{Step 2: Classification of the required orbits}.
	\vspace{0.5em}

	Using Proposition \ref{lemmadefloeser}, equation (\ref{1}) and equation (\ref{2}), we have
	\begin{align}\label{3}
	[\mathring{E_{i}}\cap \phi^{-1}(e)]=\frac{\sum_{\lambda_{1}=\dots=\lambda_{i-1}=0,\lambda_{i}=\dots=\lambda_{r}=1}[\mathcal{C}_{\lambda}^{l} \cap \mathcal{D}_{l,e}]}{(\mathbb{L}-1)\mathbb{L}^{mnl-1}} \cdot \mathbb{L}^{\nu_{i}-1}.
	\end{align}
	
	Now it suffices to calculate $[\mathcal{C}_{\lambda}^{l} \cap \mathcal{D}_{l,e}]$ for all $\lambda$ satisfying $\lambda_{1}=\dots=\lambda_{i-1}=0$ and $\lambda_{i}=\dots=\lambda_{r}=1$. 
	We consider the following map
	\begin{flalign*}
	\alpha: \mathcal{C}_{\lambda}^{l} &\rightarrow \hat{X}  \\
	             \gamma     &\mapsto  (\gamma(0),[\frac{f_{1}(\gamma)}{t^{r-i+1}}|_{t=0},\dots,\frac{f_{w}(\gamma)}{t^{r-i+1}}|_{t=0}]),
	\end{flalign*}
	where we use $t$ as the variable of the formal power series $f_{1}(\gamma),...,f_{w}(\gamma)$. This map is well-defined since the minimum order of $f_{i}(\gamma)$ is $r-i+1$. If a point $e \in \hat{X}$ lies in the image of $\alpha$, then $\alpha^{-1}(e)=\mathcal{C}_{\lambda}^{l}\cap \mathcal{D}_{l,e}$. There is a natural transitive action of $J_{l}(\mathrm{GL}_{m}) \times J_{l}(\mathrm{GL}_{n})$ on $\mathcal{C}_{\lambda}^{l}$. Moreover, this action can be naturally extended to $\alpha(\mathcal{C}_{\lambda}^{l})$ in the obvious way such that these two actions are equivariant, thus any two fibres of $\alpha$ are isomorphic under this action. This implies that
	\begin{align}\label{4}
	[\mathcal{C}_{\lambda}^{l} \cap \mathcal{D}_{l,e}]=\frac{[\mathcal{C}_{\lambda}^{l}]}{[\alpha(\mathcal{C}_{\lambda}^{l})]}.
	\end{align}
	
	\begin{remark}\label{existence_of_a_q}
	According to the proof of Proposition \ref{alsofibre}, if $\mathring{E_{i}} \cap \phi^{-1}(e) \neq \emptyset$, there exists (not necessarily unique) $\lambda$ satisfying $\lambda_{1}=\dots=\lambda_{i-1}=0,\lambda_{i}=\dots=\lambda_{r}=1$ and $e \in \alpha(\mathcal{C}_{\lambda}^{l})$.
	\end{remark}
	
	Next we want to determine which $\alpha(\mathcal{C}_{\lambda}^{l})$ contains $e$. Since the action of $J_{l}(\mathrm{GL}_{m}) \times J_{l}(\mathrm{GL}_{n})$ on $\alpha(\mathcal{C}_{\lambda}^{l})$ is transitive, for any two different $\lambda,\lambda'$, we have $\alpha(\mathcal{C}_{\lambda}^{l})=\alpha(\mathcal{C}_{\lambda'}^{l})$ or $\alpha(\mathcal{C}_{\lambda}^{l}) \cap \alpha(\mathcal{C}_{\lambda'}^{l}) =\emptyset$. For any fixed $\lambda$ such that $\lambda_{1}=\dots=\lambda_{i-1}=0$ and $\lambda_{i}=\dots=\lambda_{r}=1$, there exists a unique integer $u_{\lambda}$ such that $\lambda_{i}=\dots=\lambda_{u_{\lambda}}=1$ and $\lambda_{u_{\lambda}+1}>1$. For any $r \le q \le m$, we define
	\begin{flalign*}
		A_{q}:=\{\lambda\; |\; u_{\lambda}=q\}.
	\end{flalign*}
	We will show that these $A_{q}$ classify the different values of $\alpha(\mathcal{C}_{\lambda}^{l})$.
	\begin{proposition}\label{classification_of_c_lambda}
	For any $\lambda$ and $\lambda'$, $\alpha(\mathcal{C}_{\lambda}^{l})=\alpha(\mathcal{C}_{\lambda'}^{l})$ if and only if $\lambda$ and $\lambda'$ belong to the same $A_{q}$.
	\end{proposition}
	\begin{proof}
	If $\lambda$ and $\lambda'$ belong to the same $A_{q}$, then by the definition of $\alpha$, $\alpha(\delta_{\lambda})=\alpha(\delta_{\lambda'})$, therefore it suffices to prove the other side. Suppose $\lambda \in A_{q}$ and $\lambda' \in A_{q'}$ such that $q'>q$, then there exist $(g,h) \in J_{l}(\mathrm{GL}_{m}) \times J_{l}(\mathrm{GL}_{n})$ such that $\alpha(g\delta_{\lambda} h)=\alpha(\delta_{\lambda'})$. Noticed that the action of an element $(g,h) \in J_{l}(\mathrm{GL}_{m}) \times J_{l}(\mathrm{GL}_{n})$ on $\alpha(\mathcal{C}_{\lambda}^{l})$ is determined by their information on $0$-th jet level, therefore we have an action of $\mathrm{GL}_{m} \times \mathrm{GL}_{n}$ on $\alpha(\mathcal{C}_{\lambda}^{l})$. This action can be restated as follows.
	\begin{flalign*}
		\mathrm{GL}_{m} \times \mathrm{GL}_{n} \times \alpha(\mathcal{C}_{\lambda}^{l}) &\rightarrow \alpha(\mathcal{C}_{\lambda}^{l})    \\
		(P,Q) \cdot \alpha(\delta_{\lambda}) &\mapsto (P\begin{pmatrix}
			I_{i-1} & 0 \\
			0 & 0
		\end{pmatrix}Q,[\frac{f_{1}(G)}{t^{r-i+1}}|_{t=0},\dots,\frac{f_{w}(G)}{t^{r-i+1}}|_{t=0}]),
	\end{flalign*}
	where $G:=P\delta_{\lambda}Q$, $f_{1}(G),\dots, f_{w}(G)$ denote the value of $r \times r$ minors of $G$ and we use $t$ as the variable of the formal power series $f_{1}(G),....,f_{w}(G)$. For convenience, we will use $G^{k_{1},\dots, k_{r}}_{j_{1},\dots, j_{r}}$ to denote the $r \times r$ minor of $G$ given by $k_{1},\dots, k_{r}$-th row and $j_{1},\dots, j_{r}$-th column, where $1 \le k_{1}<\dots< k_{r} \le m$ and $1 \le j_{1}<\dots< j_{r} \le n$. The condition $(P, Q) \cdot \alpha(\delta_{\lambda})=\alpha(\delta_{\lambda'})$ is equivalent to:
	
	(i)  $P\begin{pmatrix}
		I_{i-1} & 0 \\
		0 & 0
	\end{pmatrix}Q=\begin{pmatrix}
		I_{i-1} & 0 \\
		0 & 0
	\end{pmatrix},$     
	
	(ii) For $i \le k_{i}=j_{i}<\dots<k_{r}=j_{r} \le q'$, we have $\frac{G^{1,\dots,i-1,k_{i},\dots,k_{r}}_{1,\dots,i-1,j_{i},\dots,j_{r}}}{t^{r-i+1}}|_{t=0}=a$ for some non-zero $a$, and for other minors this value is $0$. 
	
	Assume $P=\begin{pmatrix}
		P_{1} & P_{2} \\
		P_{3} & P_{4}
	\end{pmatrix}$ and  $Q=\begin{pmatrix}
		Q_{1} & Q_{2} \\
		Q_{3} & Q_{4}
	\end{pmatrix}$, then the condition (i) is equivalent to $P_{1}Q_{1}=I_{i-1}$ and $P_{3}=Q_{2}=0$. Since $\delta_{\lambda}=\begin{pmatrix}
		I_{i-1} & 0 \\
		0 & 0
	\end{pmatrix}+t \cdot \begin{pmatrix}
		0 & 0 \\
		0 & M
	\end{pmatrix}$ for some $M$, we have 
	\begin{flalign*}
		P \delta_{\lambda} Q=\begin{pmatrix}
			I_{i-1} & 0 \\
			0 & 0
		\end{pmatrix}+t \cdot \begin{pmatrix}
			P_{2}MQ_{3} & P_{2}MQ_{4} \\
			P_{4}MQ_{3} & P_{4}MQ_{4}
		\end{pmatrix}.
	\end{flalign*}
	This yields that when we calculate the $r \times r$ minors of $P \delta_{\lambda} Q$, if $\frac{G^{k_{1},\dots,k_{r}}_{j_{1},\dots,j_{r}}}{t^{r-i+1}}|_{t=0} \neq 0$, we must have $k_{1}=j_{1}=1, k_{2}=j_{2}=2,\dots,k_{i-1}=j_{i-1}=i-1$. Also, the non-zero contribution to $\frac{G^{k_{1},\dots,k_{r}}_{j_{1},\dots,j_{r}}}{t^{r-i+1}}|_{t=0}$ only comes from the constant term of $M$, therefore in this calculation we may assume $\lambda_{q+1}=\dots=\lambda_{m}=\infty$. Now the condition (ii) has been transformed into the condition for the constant matrix $G':=P_{4}MQ_{4}$ such that $G'^{k_{1},\dots,k_{r-i+1}}_{j_{1},\dots,j_{r-i+1}}=a$ if $1 \le k_{1}=j_{1}<\dots<k_{r-i+1}=j_{r-i+1} \le q'-i+1$, and otherwise $G'^{k_{1},\dots,k_{r-i+1}}_{j_{1},\dots,j_{r-i+1}}=0$. We set $G'=(b_{k,j}')_{(m-i+1) \times (n-i+1)}$, and we claim that $b_{k,j}'=0$ if $k>q'-i+1$ or $j>q'-i+1$. In fact, if $b_{1,q'-i+2}' \neq 0$, then we have
	\begin{flalign*}
		0=&det \begin{pmatrix}
			b_{1,1}' &  \cdots & b_{1,r-i+1}' & b_{1,q'-i+2}' \\
			b_{1,1}' &  \cdots & b_{1,r-i+1}' & b_{1,q'-i+2}' \\
			\vdots   &  \ddots &  \vdots     &  \vdots     \\
			b_{r-i+1,1}' &  \cdots & b_{r-i+1,r-i+1}' & b_{r-i+1,q'-i+2}' 
		\end{pmatrix}  \\
		=&\sum_{j=1}^{r-i+1}b_{1,j}'G'^{1,\dots,r-i+1}_{1,\dots,\hat{j},\dots,r-i+1,q'-i+2}(-1)^{j-1}+b_{1,q'-i+2}'G'^{1,\dots,r-i+1}_{1,\dots,r-i+1}(-1)^{r-i+1},
	\end{flalign*}
	where $\hat{j}$ means $j$ is not there. This gives a contradiction since $G'^{1,\dots,r-i+1}_{1,\dots,\hat{j},\dots,r-i+1,q'-i+2}=0$ and $G'^{1,\dots,r-i+1}_{1,\dots,r-i+1} \neq 0$. Similarly we can prove that $b_{k,j}'=0$ for $1 \le k \le q'-i+1$, $j>q'-i+1$ or $k>q'-i+1, 1 \le j \le q'-i+1$. If $b_{q'-i+2,q'-i+2} \neq 0$, then for any $1 \le k_{1}<\dots<k_{r-i} \le q'-i+1$ and $1 \le j_{1}<\dots<j_{r-i} \le q'-i+1$, $G'^{k_{1},\dots,k_{r-i}}_{j_{1},\dots,j_{r-i}}=\frac{G'^{k_{1},\dots,k_{r-i},q'-i+2}_{j_{1},\dots,j_{r-i},q'-i+2}}{b_{q'-i+2,q'-i+2}'}=0$, a contradiction. Similarly we have $b_{k,j}'=0$ for $k,j>q'-i+1$. In conclusion, $G'$ is of the form $\begin{pmatrix}
		G'' & 0 \\
		0 & 0
	\end{pmatrix}$ for some $(q'-i+1) \times (q'-i+1)$ matrix $G''$. We claim that $G''$ is a nonzero multiple of the identity matrix, which will contradict with $(P,Q) \cdot \alpha(\mathcal{C}_{\lambda}^{l})=\alpha(\delta_{\lambda'})$. In fact, since $q'>q \ge r$, we can consider the equation
	\begin{flalign*}
		0=&det \begin{pmatrix}
			b_{2,1}' &  \cdots & b_{2,r-i+2}'  \\
			b_{2,1}' &  \cdots & b_{2,r-i+2}'  \\
			b_{3,1}' &  \cdots & b_{3,r-i+2}'  \\
			\vdots   &  \ddots &  \vdots       \\
			b_{r-i+2,1}' &  \cdots & b_{r-i+2,r-i+2}'  
		\end{pmatrix}   \\
		=&\sum_{j=1}^{r-i+2}b_{2,j}'G'^{2,3,\dots,r-i+2}_{1,\dots,\hat{j},\dots,r-i+2}(-1)^{j-1}.
	\end{flalign*}
	This implies that $b_{2,1}'=0$. Similarly we can prove that $b_{k,j}'=0$ for $k \neq j$, thus $G''$ is a diagonal matrix. Now the claim follows because all the value of $(r-i+1) \times (r-i+1)$ principal minors are equal.
	\end{proof}
	
	Now for any fixed $e \in E$ such that $\mathring{E_{i}} \cap \phi^{-1}(e) \neq \emptyset$, suppose $e \in \alpha(\mathcal{C}_{\lambda}^{l})$ for some $\lambda$ and $u_{\lambda}=q$, then by Proposition \ref{classification_of_c_lambda}, equation (\ref{3}) and equation (\ref{4}), we obtain
	\begin{align}\label{5}
	[\mathring{E_{i}}\cap \phi^{-1}(e)]=\sum_{\lambda_{1}=\dots=\lambda_{i-1}=0,\; \lambda_{i}=\dots=\lambda_{r}=1,\;\lambda \in A_{q}}\frac{[\mathcal{C}_{\lambda}^{l}]}{(\mathbb{L}-1)\mathbb{L}^{mnl-\nu_{i}}[\alpha(\mathcal{C}_{\lambda}^{l})]}.
	\end{align}

	\vspace{0.5em}
	\noindent\textbf{Step 3: Calculation of $[\mathcal{C}_{\lambda}^{l}]$ and $[\alpha(\mathcal{C}_{\lambda}^{l})]$}.
	\vspace{0.5em}
	
	We will finish our proof by computing the value of $[\mathcal{C}_{\lambda}^{l}]$ and $[\alpha(\mathcal{C}_{\lambda}^{l})]$ introduced above. Before that, we need some preparations.
	\begin{definition}
		Assume $0 < v_{1}<\dots<v_{j}<m$ are fixed integers. We call a chain of $k$-vector spaces $V_{1} \subset V_{2} \subset \dots \subset V_{j} \subset \mathbb{C}^{m}$ a flag if the dimension of $V_{i}$ is $v_{i}$, and $(v_{1},\dots,v_{j})$ is called the signature of the flag. Note that if the signature is fixed, then $\mathrm{GL}_{m}$ acts transitively on the set of flags, and we call the stabilizer of a flag parabolic subgroup of $\mathrm{GL}_{m}$. If $P$ is a parabolic subgroup, then $\mathrm{GL}_{m}/P$ parametrizes the flag of a given signature.
	\end{definition}
	\begin{definition}\label{Definition of P and L}
		Suppose $\{e_{1},\dots,e_{m}\}$ is a basis for $\mathbb{C}^{m}$, and $$\lambda=(d_{1},\dots,d_{1},d_{2},\dots,d_{2},\dots,d_{j+1},\dots,d_{j+1})$$ is a vector with $c_{i}$ many $d_{i}$, $0 \le d_{1}<\dots<d_{j+1}$ and $c_{1}+\dots+c_{j+1}=m$. We set $v_{i}=c_{1}+\dots+c_{i}$ for $1 \le i \le j+1$, and $V_{i}=\mathrm{Span}\{e_{1},\dots,e_{v_{i}}\}$, then $V_{1} \subset\dots \subset V_{j} \subset V_{j+1}=\mathbb{C}^{m}$ is a flag. We denote the stabilizer of this flag $P_{\lambda}$ and we call it the parabolic subgroup associated to $\lambda$. $L_{\lambda}:=\mathrm{GL}_{c_{1}} \times \mathrm{GL}_{c_{2}} \times \dots \times \mathrm{GL}_{c_{j+1}}$ is a subgroup of $P_{\lambda}$, and we call it the Levi factor associated to $\lambda$.
	\end{definition}
	
	If we fix a basis $\{e_{1},\dots,e_{m}\}$ as in the definition above and express $P_{\lambda}, L_{\lambda}$ in matrices with respect to this basis, $P_{\lambda}$ will be those upper triangular matrices in blocks and $L_{\lambda}$ will be those diagonal matrices in blocks. The following example describes this explicitly.
	\begin{example}
		Suppose $m=5$, $\lambda=(0,1,1,3,3)$, $c_{1}=1,c_{2}=2,c_{3}=2$, then $v_{1}=1,v_{2}=3,v_{3}=5$. $P_{\lambda}$ and $L_{\lambda}$ consist of the matrices of the form
		\begin{flalign*}
			P_{\lambda}=\begin{pmatrix}
				* & * & * & * & * \\
				0 & * & * & * & * \\
				0 & * & * & * & * \\
				0 & 0 & 0 & * & * \\
				0 & 0 & 0 & * & * 
			\end{pmatrix}, \; \;L_{\lambda}=\begin{pmatrix}
				* & 0 & 0 & 0 & 0 \\
				0 & * & * & 0 & 0 \\
				0 & * & * & 0 & 0 \\
				0 & 0 & 0 & * & * \\
				0 & 0 & 0 & * & * 
			\end{pmatrix}.
		\end{flalign*}
	\end{example}
	
	Now we give the value of $[\mathcal{C}_{\lambda}^{l}]$.
	\begin{proposition}\label{c_lambda}
	We assume $\lambda=(\lambda_{1},\dots,\lambda_{m})$ with $\lambda_{1}=\dots=\lambda_{i-1}=0$ and $\lambda_{i}=\dots=\lambda_{r}=1$. Suppose $\lambda_{r+1} \le \dots \le \lambda_{s} \le l$ and $\lambda_{s+1}=\dots=\lambda_{m}=\infty$, where $r \le s \le m$. Let \begin{flalign*}
		\lambda_{1}&=a_{1}, \\
		\lambda_{2}&=a_{1}+a_{2}, \\
		&\vdots \\
		\lambda_{s}&=a_{1}+\dots+a_{s},
	\end{flalign*}
	and we define $I:=\{2 \le i \le s \; | \; a_{i} \neq 0\}$. Suppose $I \cup \{1\}=\{i_{1},\dots,i_{p}\}$, where $i_{1}>\dots>i_{p}$, then we have
	\begin{flalign*}
		[\mathcal{C}_{\lambda}^{l}]=\frac{[\mathrm{GL}_{m}] \cdot [\mathrm{GL}_{n}] \cdot \prod_{j=1}^{p}[G(i_{j-1}-i_{j},s+1-i_{j})]^{2} \cdot [\mathrm{GL}_{i_{j-1}-i_{j}}]}{[\mathrm{GL}_{s}]^{2} \cdot [\mathrm{GL}_{m-s}] \cdot [\mathrm{GL}_{n-s}] \cdot \mathbb{L}^{s^{2}l-(m+n)sl+s(m+n-2s)+\sum_{j=1}^{s}(m+n+1-2j)\lambda_{j}}},
	\end{flalign*}
	where we define $i_{0}=s+1$. 
	\end{proposition}
	
	\begin{proof}
	 Since we have the action of $J_{l}(\mathrm{GL}_{m}) \times J_{l}(\mathrm{GL}_{n})$ on $\mathcal{C}_{\lambda}^{l}$, we have 
	\begin{flalign*}
		[\mathcal{C}_{\lambda}^{l}]=\frac{[J_{l}(\mathrm{GL}_{m}) \times J_{l}(\mathrm{GL}_{n})]}{[H_{\lambda}]},
	\end{flalign*}
	where $H_{\lambda}$ is the stabilizer of $\delta_{\lambda}$ under the action of $J_{l}(\mathrm{GL}_{m}) \times J_{l}(\mathrm{GL}_{n})$. Assume $(g_{i,j})_{m \times m} \times (h_{i,j})_{n \times n} \in H_{\lambda}$, and we set $g_{i,j}=g_{i,j}^{(0)}+g_{i,j}^{(1)}t+\dots+g_{i,j}^{(l)}t^{l},h_{i,j}=h_{i,j}^{(0)}+h_{i,j}^{(1)}t+\dots+h_{i,j}^{(l)}t^{l}$. Since $(g_{i,j}) \cdot \delta_{\lambda}=\delta_{\lambda} \cdot (h_{i,j})$, we have
	\begin{flalign*}
		\begin{pmatrix}
			t^{\lambda_{1}}g_{1,1} & \cdots & t^{\lambda_{s}}g_{1,s} & 0 & \cdots & 0 \\
			\vdots & \ddots & \vdots & \vdots & \vdots & \vdots \\
			t^{\lambda_{1}}g_{s,1} & \cdots & t^{\lambda_{s}}g_{s,s} & 0 & \cdots & 0 \\
			\vdots & \ddots & \vdots & \vdots & \vdots & \vdots \\
			t^{\lambda_{1}}g_{m,1} & \cdots & t^{\lambda_{s}}g_{m,s} & 0 & \cdots & 0
		\end{pmatrix}=
		\begin{pmatrix}
			t^{\lambda_{1}}h_{1,1} & \cdots & t^{\lambda_{1}}h_{1,s} & \cdots & t^{\lambda_{1}}h_{1,n}  \\
			\vdots & \ddots & \vdots & \ddots & \vdots \\
			t^{\lambda_{s}}h_{s,1} & \cdots & t^{\lambda_{s}}h_{s,s} & \cdots & t^{\lambda_{s}}h_{s,n} \\
			0 & \cdots & 0 & \cdots & 0 \\
			\vdots & \vdots & \vdots & \vdots & \vdots \\
			0 & \cdots & 0 & \cdots & 0
		\end{pmatrix}.
	\end{flalign*}
	This implies that $t^{\lambda_{j}}g_{i,j}=0$ for $s+1 \le i \le m, 1 \le j \le s$ and $t^{\lambda_{i}}h_{i,j}=0$ for $1 \le i \le s, s+1  \le j \le n$, therefore $g_{i,j}^{(0)}=0$ and $h_{i,j}^{(0)}=0$ for these corresponding $i,j$. Since the condition $(g_{i,j})_{m \times m} \times (h_{i,j})_{n \times n} \in J_{l}(\mathrm{GL}_{m}) \times J_{l}(\mathrm{GL}_{n})$ is equivalent to $(g_{i,j}^{(0)})_{m \times m} \times (h_{i,j}^{(0)})_{n \times n} \in \mathrm{GL}_{m} \times \mathrm{GL}_{n}$, this restriction is equivalent to 
	\begin{flalign*}
		&(g_{i,j}^{(0)}\;|\;1 \le i,j \le s) \in \mathrm{GL}_{s},  \\ &(g_{i,j}^{(0)}\;| \; s+1 \le i,j \le m) \in \mathrm{GL}_{m-s},  \\
		&(h_{i,j}^{(0)}\;|\;1 \le i,j \le s) \in \mathrm{GL}_{s}, \\
		&(h_{i,j}^{(0)}\;|\; s+1 \le i,j \le n) \in \mathrm{GL}_{n-s}.
	\end{flalign*}
	Now if we define $H_{\lambda}'$ to be the set of pairs $(g_{i,j}\;|\;1 \le i,j \le s), (h_{i,j}\;|\;1 \le i,j \le s)$ such that
	\begin{flalign*}
		\begin{pmatrix}
			t^{\lambda_{1}}g_{1,1} & \cdots & t^{\lambda_{s}}g_{1,s}  \\
			\vdots & \ddots & \vdots  \\
			t^{\lambda_{1}}g_{s,1} & \cdots & t^{\lambda_{s}}g_{s,s}  \\
		\end{pmatrix}=
		\begin{pmatrix}
			t^{\lambda_{1}}h_{1,1} & \cdots & t^{\lambda_{1}}h_{1,s} \\
			\vdots & \ddots & \vdots \\
			t^{\lambda_{s}}h_{s,1} & \cdots & t^{\lambda_{s}}h_{s,s} \\
		\end{pmatrix},
	\end{flalign*}
	then we have
	\begin{flalign*}
		[H_{\lambda}]=&[J_{l}(\mathrm{GL}_{m-s})] \cdot [J_{l}(\mathrm{GL}_{n-s})] \cdot \mathbb{L}^{(m+n-2s)\sum_{j=1}^{s}\lambda_{j}+s(m-s)(l+1)+s(n-s)(l+1)} \cdot [H_{\lambda}'] \\
		=&[H_{\lambda}'] \cdot [\mathrm{GL}_{m-s}] \cdot [\mathrm{GL}_{n-s}] \cdot \mathbb{L}^{m(m-s)l+n(n-s)l+s(m+n-2s)+(m+n-2s)\sum_{j=1}^{s}\lambda_{j}}.
	\end{flalign*}
	According to Proposition $6.4$ in \cite{Arcs_on_Determinantal_Variteties_Roi}, if we define $\lambda':=(\lambda_{1},\dots,\lambda_{s})$, then we have $[H_{\lambda}']=\frac{[P_{\lambda'}]^{2} \cdot \mathbb{L}^{s^{2}l+\sum_{j=1}^{s}(2s+1-2j)\lambda_{j}}}{[L_{\lambda'}]}$. 
	By the calculation of Section $6.2$ in \cite{Arcs_on_Determinantal_Variteties_Roi}, we obtain
	\begin{flalign*}
		[\mathrm{GL}_{s}/P_{\lambda'}]^{2}\cdot [L_{\lambda'}]=\prod_{j=1}^{p}[G(i_{j-1}-i_{j},s+1-i_{j})]^{2}\cdot [\mathrm{GL}_{i_{j-1}-i_{j}}],
	\end{flalign*} 
	where $i_{0}=s+1$ and $G(d,k)$ denotes the Grassmannian of $d$ dimensional subspace of $k$ dimensional vector space. This implies that
	\begin{flalign*}
		[H_{\lambda}]=\frac{[\mathrm{GL}_{s}]^{2} \cdot [\mathrm{GL}_{m-s}] \cdot [\mathrm{GL}_{n-s}] \cdot \mathbb{L}^{s^{2}l+m(m-s)l+n(n-s)l+s(m+n-2s)+\sum_{j=1}^{s}(m+n+1-2j)\lambda_{j}}}{\prod_{j=1}^{p}[G(i_{j-1}-i_{j},s+1-i_{j})]^{2} \cdot [\mathrm{GL}_{i_{j-1}-i_{j}}]}
	\end{flalign*}
	and
	\begin{flalign*}
		[\mathcal{C}_{\lambda}^{l}]=\frac{[\mathrm{GL}_{m}] \cdot [\mathrm{GL}_{n}] \cdot \prod_{j=1}^{p}[G(i_{j-1}-i_{j},s+1-i_{j})]^{2} \cdot [\mathrm{GL}_{i_{j-1}-i_{j}}]}{[\mathrm{GL}_{s}]^{2} \cdot [\mathrm{GL}_{m-s}] \cdot [\mathrm{GL}_{n-s}] \cdot \mathbb{L}^{s^{2}l-(m+n)sl+s(m+n-2s)+\sum_{j=1}^{s}(m+n+1-2j)\lambda_{j}}}.
	\end{flalign*} 
	\end{proof}
	
	\vspace{-0.5 em} 
	 Next we compute $[\alpha(\mathcal{C}_{\lambda}^{l})]$. 
	 \begin{proposition}\label{alpha_c_lambda}
	 Let $q$ be the integer such that $r \le q \le s$, $\lambda_{i}=\dots=\lambda_{r}=\dots=\lambda_{q}=1$ and $\lambda_{q+1}>1$. If $q=r$, then
	 \begin{flalign*}
	 	[\alpha(\mathcal{C}_{\lambda}^{l})]=\frac{[\mathrm{GL}_{m}] \cdot [\mathrm{GL}_{n}]}{[\mathrm{GL}_{i-1}] \cdot [\mathrm{GL}_{r-i+1}]^{2} \cdot [\mathrm{GL}_{m-r}] \cdot [\mathrm{GL}_{n-r}] \cdot \mathbb{L}^{(m+n-2i+2)(i-1)+(r-i+1)(m+n-2r)}}.
	 \end{flalign*}
	 If $q>r$, then
	 \begin{flalign*}
	 	[\alpha(\mathcal{C}_{\lambda}^{l})]=\frac{[\mathrm{GL}_{m}] \cdot [\mathrm{GL}_{n}]}{[\mathrm{GL}_{i-1}] \cdot [\mathrm{GL}_{q-i+1}] \cdot [\mathrm{GL}_{m-q}] \cdot [\mathrm{GL}_{n-q}] \cdot \mathbb{L}^{(m+n-2i+2)(i-1)+(q-i+1)(m+n-2q)} \cdot (\mathbb{L}-1)}.
	 \end{flalign*}
	 \end{proposition}
	 \begin{proof}
	As in the proof of Proposition \ref{classification_of_c_lambda}, we consider the action of $\mathrm{GL}_{m} \times \mathrm{GL}_{n}$ on $\alpha(\mathcal{C}_{\lambda}^{l})$ as follows.
	 \begin{flalign*}
	 	\mathrm{GL}_{m} \times \mathrm{GL}_{n} \times \alpha(\mathcal{C}_{\lambda}^{l}) &\rightarrow \alpha(\mathcal{C}_{\lambda}^{l})    \\
	 	(P,Q) \cdot \alpha(\delta_{\lambda}) &\mapsto (P\begin{pmatrix}
	 		I_{i-1} & 0 \\
	 		0 & 0
	 	\end{pmatrix}Q,[\frac{f_{1}(G)}{t^{r-i+1}}|_{t=0},\dots,\frac{f_{w}(G)}{t^{r-i+1}}|_{t=0}]),
	 \end{flalign*}
	 where $G:=P\delta_{\lambda}Q$, $f_{1}(G),\dots,f_{w}(G)$ denote the value of $r \times r$ minors of $G$ and we use $t$ as the variable of the formal power series $f_{1}(G),...,f_{w}(G)$. We define $Z_{\lambda}$ to be the stabilizer of $\delta_{\lambda}$ in $\mathrm{GL}_{m} \times \mathrm{GL}_{n}$, then we have
	 \begin{flalign*}
	 	[\alpha(\mathcal{C}_{\lambda}^{l})]=\frac{[\mathrm{GL}_{m}] \cdot [\mathrm{GL}_{n}]}{[Z_{\lambda}]}.
	 \end{flalign*}
	 Suppose $(P,Q) \in Z_{\lambda}$, then $(P,Q) \cdot \alpha(\delta_{\lambda})=\alpha(\delta_{\lambda})$, which is equivalent to
	 
	 (i)  $P\begin{pmatrix}
	 	I_{i-1} & 0 \\
	 	0 & 0
	 \end{pmatrix}Q=\begin{pmatrix}
	 	I_{i-1} & 0 \\
	 	0 & 0
	 \end{pmatrix},$     
	 
	 (ii) For $i \le k_{i}=j_{i}<\dots<k_{r}=j_{r} \le q$, we have $\frac{G^{1,\dots,i-1,k_{i},\dots,k_{r}}_{1,\dots,i-1,j_{i},\dots,j_{r}}}{t^{r-i+1}}|_{t=0}=a$ for some non-zero $a$, and for other minors this value is $0$. 
	 
	 Again we assume $P=\begin{pmatrix}
	 	P_{1} & P_{2} \\
	 	P_{3} & P_{4}
	 \end{pmatrix}$ and  $Q=\begin{pmatrix}
	 	Q_{1} & Q_{2} \\
	 	Q_{3} & Q_{4}
	 \end{pmatrix}$, then the condition (i) is equivalent to $P_{1}Q_{1}=I_{i-1}$ and $P_{3}=Q_{2}=0$. Since $\delta_{\lambda}=\begin{pmatrix}
	 	I_{i-1} & 0 \\
	 	0 & 0
	 \end{pmatrix}+t \cdot \begin{pmatrix}
	 	0 & 0 \\
	 	0 & M
	 \end{pmatrix}$ for some $M$, we have 
	 \begin{flalign*}
	 	P \delta_{\lambda} Q=\begin{pmatrix}
	 		I_{i-1} & 0 \\
	 		0 & 0
	 	\end{pmatrix}+t \cdot \begin{pmatrix}
	 		P_{2}MQ_{3} & P_{2}MQ_{4} \\
	 		P_{4}MQ_{3} & P_{4}MQ_{4}
	 	\end{pmatrix}.
	 \end{flalign*}
	 As in the proof of Proposition \ref{classification_of_c_lambda}, the condition (ii) is equivalent to condition for the constant matrix $G':=P_{4}MQ_{4}$ such that $G'^{k_{1},\dots,k_{r-i+1}}_{j_{1},\dots,j_{r-i+1}}=a$ if $1 \le k_{1}=j_{1}<\dots<k_{r-i+1}=j_{r-i+1} \le q-i+1$, and otherwise $G'^{k_{1},\dots,k_{r-i+1}}_{j_{1},\dots,j_{r-i+1}}=0$. Again we can prove that $G'$ is of the form $\begin{pmatrix}
	 	G'' & 0 \\
	 	0 & 0
	 \end{pmatrix}$ for some $(q-i+1) \times (q-i+1)$ matrix $G''$. If $q=r$, then any invertible matrix $G''$ will satisfy the condition. If $q>r$, $G''$ will be a nonzero multiple of the identity matrix. 
	 
	 Now we can compute the value of $[Z_{\lambda}]$. If $q=r$, the restriction on $P_{4},Q_{4}$ is
	 \begin{flalign*}
	 	P_{4}\begin{pmatrix}
	 		I_{r-i+1} & 0 \\
	 		0 & 0
	 	\end{pmatrix}Q_{4}=\begin{pmatrix}
	 		G'' & 0 \\
	 		0 & 0
	 	\end{pmatrix}
	 \end{flalign*}
	 for some invertible matrix $G''$. If we write $P_{4}=\begin{pmatrix}
	 	P_{5} & P_{6} \\
	 	P_{7} & P_{8}
	 \end{pmatrix}$ and $Q_{4}=\begin{pmatrix}
	 	Q_{5} & Q_{6} \\
	 	Q_{7} & Q_{8}
	 \end{pmatrix}$, then this condition is equivalent to $P_{3}=Q_{6}=0$. In conclusion, we have
	 \begin{flalign*}
	 	[Z_{\lambda}]=[\mathrm{GL}_{i-1}] \cdot [\mathrm{GL}_{r-i+1}]^{2} \cdot [\mathrm{GL}_{m-r}] \cdot [\mathrm{GL}_{n-r}] \cdot \mathbb{L}^{(m+n-2i+2)(i-1)+(r-i+1)(m+n-2r)}
	 \end{flalign*}
	 and
	 \begin{flalign*}
	 	[\alpha(\mathcal{C}_{\lambda}^{l})]=\frac{[\mathrm{GL}_{m}] \cdot [\mathrm{GL}_{n}]}{[\mathrm{GL}_{i-1}] \cdot [\mathrm{GL}_{r-i+1}]^{2} \cdot [\mathrm{GL}_{m-r}] \cdot [\mathrm{GL}_{n-r}] \cdot \mathbb{L}^{(m+n-2i+2)(i-1)+(r-i+1)(m+n-2r)}}.
	 \end{flalign*}
	 
	 If $q>r$, the restriction on $P_{4},Q_{4}$ is
	 \begin{flalign*}
	 	P_{4}\begin{pmatrix}
	 		I_{q-i+1} & 0 \\
	 		0 & 0
	 	\end{pmatrix}Q_{4}=\begin{pmatrix}
	 		aI_{q-i+1} & 0 \\
	 		0 & 0
	 	\end{pmatrix}
	 \end{flalign*}
	 for some nonzero $a$. Again if we write $P_{4}=\begin{pmatrix}
	 	P_{5} & P_{6} \\
	 	P_{7} & P_{8}
	 \end{pmatrix}$ and $Q_{4}=\begin{pmatrix}
	 	Q_{5} & Q_{6} \\
	 	Q_{7} & Q_{8}
	 \end{pmatrix}$, then this condition is equivalent to $P_{3}=Q_{6}=0$ and $P_{5}Q_{5}=aI_{q-i+1}$.
	 In conclusion, we have
	 \begin{flalign*}
	 	[Z_{\lambda}]=[\mathrm{GL}_{i-1}] \cdot [\mathrm{GL}_{q-i+1}] \cdot [\mathrm{GL}_{m-q}] \cdot [\mathrm{GL}_{n-q}] \cdot \mathbb{L}^{(m+n-2i+2)(i-1)+(q-i+1)(m+n-2q)} \cdot (\mathbb{L}-1)
	 \end{flalign*}
	 and
	 \begin{flalign*}
	 	[\alpha(\mathcal{C}_{\lambda}^{l})]=\frac{[\mathrm{GL}_{m}] \cdot [\mathrm{GL}_{n}]}{[\mathrm{GL}_{i-1}] \cdot [\mathrm{GL}_{q-i+1}] \cdot [\mathrm{GL}_{m-q}] \cdot [\mathrm{GL}_{n-q}] \cdot \mathbb{L}^{(m+n-2i+2)(i-1)+(q-i+1)(m+n-2q)} \cdot (\mathbb{L}-1)}.
	 \end{flalign*}
	 \end{proof}

	There are some other related varieties and their values in the Grothendieck ring are given in the following lemma.
	\begin{lemma}[Example $2.4.4$ and $2.4.5$ in \cite{Motivic_Integration_book}]\label{Grassmannian}
		For integer $0<d \le k$, we have:
		
		(i)The general linear group $[\mathrm{GL}_{d}]=\mathbb{L}^{d(d-1)/2}(\mathbb{L}^{d}-1)(\mathbb{L}^{d-1}-1) \cdots (\mathbb{L}-1)$.
		
		(ii)The Grassmannian $[G(d,k)]=\prod_{j=1}^{d}\frac{\mathbb{L}^{j+k-d}-1}{\mathbb{L}^{j}-1}=\sum_{0 \le \lambda_{1} \le \dots \le \lambda_{d} \le k-d}\mathbb{L}^{\lambda_{1}+\dots+\lambda_{d}}$.
	\end{lemma}
	
	Now we can give the monodromy zeta function of $Z_{r}$ at any point $e \in E$. This is an explicit version of Theorem B.
	
	\begin{theorem}\label{thmB'}
	For $e \in E$, $h(e) \in Z_{r}$, and we suppose the rank of $h(e)$ is $i-1$ and $1 \le i \le r$. If there is no $\lambda$ satisfying $\lambda_{1}=\dots=\lambda_{i-1}=0$ and $\lambda_{i}=\dots=\lambda_{r}=1$ such that $e \in \alpha(\mathcal{C}_{\lambda}^{l})$, then $Z^{\mathrm{mon}}_{Z_{r},e}=1$. Otherwise $e \in \alpha(\mathcal{C}_{\lambda}^{l})$ for some $\lambda \in A_{q}$. If $i<r=q$, then $Z^{\mathrm{mon}}_{Z_{r},e}$ is still $1$. If $i=r$ or $i<r<q$, then $Z^{\mathrm{mon}}_{Z_{r},e}=1-t^{r+1-i}$. 
	\end{theorem}
	\begin{proof}
	By Remark \ref{only_one_e_nonzero}, only when $j=i$ $\mathring{E_{i}} \cap \phi^{-1}(e)$ may not be empty, therefore $Z^{\mathrm{mon}}_{Z_{r},e}=(1-t^{N_{i}})^{-\chi(\mathring{E_{i}} \cap \phi^{-1}(e))}$. If $\mathring{E_{i}} \cap \phi^{-1}(e)$ is not empty, by Remark \ref{existence_of_a_q} we can suppose $e \in \alpha(\mathcal{C}_{\lambda}^{l})$ and $\lambda \in A_{q}$. By Lemma \ref{Grassmannian}, the power of $\mathbb{L}-1$ in $[\mathrm{GL}_{d}],G(d,k)$ is $d$ and $0$, respectively. This implies that the power of $\mathbb{L}-1$ in $\mathcal{C}_{\lambda}^{l}$ is $s$ and that of $[\alpha(\mathcal{C}_{\lambda}^{l})]$ is $i-1$ or $q-1$. Now we analyze the power of every term in the equation (\ref{5}). Since $s \ge q \ge r \ge i$, the power of $\mathbb{L}-1$ in $\frac{\mathcal{C}_{\lambda}^{l}}{\alpha(\mathcal{C}_{\lambda}^{l})(\mathbb{L}-1)}$ is zero if and only if $s=q=r=i$ or $s=q>r \ge i$. In both cases, there is only one $\lambda$ in $A_{q}$ satisfying these conditions. When $s=q>r \ge i$, direct computation with the value given in Proposition \ref{c_lambda}, Proposition \ref{alpha_c_lambda} and Lemma \ref{Grassmannian} shows
	\begin{flalign*}
		&\chi(\frac{\mathcal{C}_{\lambda}^{l}}{\alpha(\mathcal{C}_{\lambda}^{l})(\mathbb{L}-1)})   \\
		=&\frac{m! \cdot n! \cdot (q!)^{2}}{(q!)^{2} \cdot (m-q)! \cdot (n-q)! \cdot (i-1)! \cdot (q-i+1)!}  \cdot \frac{(i-1)! \cdot (q-i+1)! \cdot (m-q)! \cdot (n-q)!}{m! \cdot n!}  \\
		=&1.
	\end{flalign*}
	Similarly when $s=q=r=i$, we have
	\begin{flalign*}
		&\chi(\frac{\mathcal{C}_{\lambda}^{l}}{\alpha(\mathcal{C}_{\lambda}^{l})(\mathbb{L}-1)})   \\
		=&\frac{m! \cdot n! \cdot (q!)^{2}}{(q!)^{2} \cdot (m-q)! \cdot (n-q)! \cdot (q-1)!}  \cdot \frac{(q-1)! \cdot (m-q)! \cdot (n-q)!}{m! \cdot n!}  \\
		=&1.
	\end{flalign*}
	In conclusion, only if $i=r$ or $i<r<q$, there is one $\lambda$ in the equation (\ref{5}) such that $\chi(\frac{\mathcal{C}_{\lambda}^{l}}{\alpha(\mathcal{C}_{\lambda}^{l})(\mathbb{L}-1)})$ is non-zero and $Z^{\mathrm{mon}}_{Z_{r},e}=1-t^{r+1-i}$. In other cases $Z^{\mathrm{mon}}_{Z_{r},e}=1$. 
	\end{proof}

	\begin{theorem}\label{thmB}
	Let $\hat{X} \rightarrow X$ be the blow-up of $X$ at $Z_{r}$ and let $E$ be its exceptional divisor. For different points $e \in E$, the monodromy zeta function of $Z_{r}$ at $e$, denoted by $Z^{\mathrm{mon}}_{Z_{r},e}(t)$, can be $1$ or $1-t^{r+1-i}$ for $1 \le i \le r$, depending on the point $e$.
	\end{theorem}
	\begin{proof}
	This follows directly from Theorem \ref{thmB'}.
	\end{proof}

	\begin{theorem}\label{thmA}
	The monodromy conjecture (Conjecture \ref{mc}) holds for $Z_{r}$.
	\end{theorem}
	\begin{proof}
	According to Remark \ref{poles_of_zeta_function}, the poles of the motivic zeta function of $Z_{r}$ are given by $-\frac{(m+1-j)(n+1-j)}{r+1-j}$. On the other hand, the zeros and poles of the monodromy zeta function, which by Proposition \ref{ideal_case_monodromy_zeta_function_and_eigenvalue} can fully characterize the eigenvalues of the Verdier monodromy, is given in Theorem \ref{thmB}. The theorem follows by comparing these two results. 
	\end{proof}
	 
	\begin{example}
	Let us see the simplest case where monodromy zeta function may have nontrivial zeros and poles. Let $m=n=3$, $r=2$, and $X=\mathbb{C}^{9}=\mathrm{Spec}\mathbb{C}[x_{11},\dots,x_{33}]$. The log resolution of $(X,Z_{2})$ consists of two blow-ups. Firstly we blow up $0$ and obtain $\pi_{1}:X_{1} \rightarrow X$, where $X_{1}=\{((x_{11},\dots,x_{33}),[u_{11},\dots,u_{33}]) \in \mathbb{C}^{9} \times \mathbb{P}^{8}\;|\; x_{ij}=\eta u_{ij}, \forall i,j\}$. Suppose $U_{ij}=\{u_{ij} \neq 0\}$ is a chart of $X_{1}$ and all the $U_{ij}$ form an open covering of $X_{1}$. Let us focus on the chart $U_{11}$. On this chart we may assume $u_{11}=0$ and choose $\eta,u_{12},u_{13},u_{21},\dots,u_{33}$ as coordinates. The strict transform $\tilde{Z}_{2}$ is defined by
	\begin{flalign*}
		u_{22}-u_{12}u_{21},u_{23}-u_{13}u_{21},u_{32}-u_{12}u_{31},u_{33}-u_{13}u_{31}.
	\end{flalign*}
	Let $u_{ij}':=u_{ij}-u_{1j}u_{i1}$ for $2 \le i,j \le 3$, then these $u_{ij}'$ together with $\eta,u_{12},u_{13},u_{21},u_{31}$ form a new set of coordinates. Next we blow up $\tilde{Z}_{2}$ and obtain $\pi_{2}: X_{2} \rightarrow X_{1}$, where $\tilde{U}_{11}:=X_{2} \cap \pi_{2}^{-1}(U_{11})=\{(\eta,u_{12},u_{13},u_{21},u_{31},u_{22}',u_{23}',u_{32}',u_{33}'),[v_{22}',v_{23}',v_{32}',v_{33}']\; | \; u_{ij}'=\theta v_{ij}'\}$.  In this chart $E_{1},E_{2}$ are defined by $\eta, \theta$, respectively. Let $h: \hat{X} \rightarrow X$ be the blow-up of $Z_{2}$, $E$ be its exceptional divisor and $\phi: X_{2} \rightarrow \hat{X}$ be the canonical morphism such that $h \circ \phi=\pi_{1} \circ \pi_{2}$. 
	
	Next we will calculate $Z^{mon}_{Z_{2},e}(t)$ for all possible $e \in E$. By Proposition \ref{classification_of_c_lambda}, to compute possible non-zero $\chi(\mathring{E_{i}} \cap \phi^{-1}(e))$, we only need to consider $e=\alpha(\delta_{\lambda})$ for some $\lambda$ satisfying $\lambda_{1}=\dots=\lambda_{i-1}=0,\lambda_{i}=\dots\lambda_{2}=1$. 
	
	When $i=1$, it suffices to consider $e=\alpha(\delta_{\lambda})$ where $\lambda=(1,1,1)$ or $\lambda=(1,1,\infty)$. When $\lambda=(1,1,\infty)$, $e \in  E$ is defined by $x_{ij}=0$ and $[f_{1}(x_{ij}),\dots,f_{9}(x_{ij})]=[1,0,\dots,0]$, where $f_{1},\dots,f_{9}$ denote the $2 \times 2$ minors of the matrix $(x_{ij})_{3 \times 3}$ and $f_{1}=x_{11}x_{22}-x_{12}x_{21}$, i.e. $e$ is in the chart $\mathrm{Spec}\mathbb{C}[x_{11},\dots,x_{33},\frac{f_{2}}{f_{1}},\dots,\frac{f_{9}}{f_{1}}]$ and is given by $0$. Direct computation shows that $\mathring{E_{1}} \cap \phi^{-1}(e) \cap \tilde{U}_{11}$ is given by
	\begin{flalign*}
	\eta=0, \theta \neq 0, [v_{22}',v_{23}',v_{32}',v_{33}']=[1,0,0,0],u_{13}=u_{31}=0.
	\end{flalign*}
	This yields $[\mathring{E_{1}} \cap \phi^{-1}(e) \cap \tilde{U}_{11}]=(\mathbb{L}-1)\mathbb{L}^{2}$. Similarly we can consider the chart $\tilde{U}_{12}:=X_{2} \cap \pi_{2}^{-1}(U_{12})=\{(\eta,u_{11},u_{13},u_{22},u_{32},u_{21}',u_{23}',u_{31}',u_{33}'),[v_{21}',v_{23}',v_{31}',v_{33}']\; | \; u_{ij}'=\theta v_{ij}'\}$. Direct calculation shows $\mathring{E_{1}} \cap \phi^{-1}(e) \cap (\tilde{U}_{12}-\tilde{U}_{11})$ is given by
	\begin{flalign*}
		\eta=0, \theta \neq 0, [v_{21}',v_{23}',v_{31}',v_{33}']=[1,0,0,0], u_{11}=u_{13}=u_{32}=0,
	\end{flalign*}
	therefore we have $[\mathring{E_{1}} \cap \phi^{-1}(e) \cap (\tilde{U}_{12}-\tilde{U}_{11})]=(\mathbb{L}-1)\mathbb{L}$. Similarly we can consider other charts and finally we can prove that $\mathring{E_{1}} \cap \phi^{-1}(e)$ is covered by $\tilde{U}_{11}$ and $\tilde{U}_{12}$, thus $\chi(\mathring{E_{1}} \cap \phi^{-1}(e))=0$. This yields $Z^{\mathrm{mon}}_{Z_{2},e}=1$.
	
	When $\lambda=(1,1,1)$, $e \in E$ is defined by $x_{ij}=0$ and $f_{1}=f_{2}=f_{3} \neq 0,f_{4}=\dots=f_{9}=0$, where $f_{1}=x_{11}x_{22}-x_{12}x_{21},f_{2}=x_{11}x_{33}-x_{13}x_{31},f_{3}=x_{22}x_{33}-x_{23}x_{32}$. Direct computation shows that $\mathring{E_{1}} \cap \phi^{-1}(e) \cap \tilde{U}_{11}$ is given by
	\begin{flalign*}
	\eta=0, u_{12}=u_{13}=u_{21}=u_{31}=u_{23}'=u_{32}'=0,u_{22}'=u_{33}'=1,
	\end{flalign*}
	i.e. it is a point. Similarly we can see that $\mathring{E_{1}} \cap \phi^{-1}(e)$ is covered by $\tilde{U}_{11}$, thus $\chi(\mathring{E_{1}} \cap \phi^{-1}(e))=1$ and $Z^{\mathrm{mon}}_{Z_{2},e}=1-t^{2}$. 
	
	Finally we study the case when $i=2$. Similarly it suffices to consider $e=\alpha(\delta_{\lambda})$, where $\lambda=(0,1,1)$ or $\lambda=(0,1,\infty)$. When $\lambda=(0,1,\infty)$, $e \in E$ is determined by $x_{11}=1$, $x_{ij}=0$ for $(i,j) \neq (1,1)$, and $[f_{1},\dots,f_{9}]=[1,0,\dots,0]$. This implies that $\mathring{E_{2}} \cap \phi^{-1}(e) \cap \tilde{U}_{11}$ is given by
	\begin{flalign*}
	\eta=1,u_{12}=u_{13}=u_{21}=u_{31}=\theta=v_{23}'=v_{32}'=v_{33}'=0,v_{22}'=1,
	\end{flalign*}
	i.e. $\mathring{E_{1}} \cap \phi^{-1}(e) \cap \tilde{U}_{11}$ is a point. Since $x_{11}=1$, $\mathring{E_{1}} \cap \phi^{-1}(e) \subset \tilde{U}_{11}$, thus $\chi(\mathring{E_{1}} \cap \phi^{-1}(e))=1$ and $Z^{\mathrm{mon}}_{Z_{2},e}=1-t$.
	
	When $\lambda=(0,1,1)$, $e \in E$ is determined by $x_{11}=1$, $x_{ij}=0$ for $(i,j) \neq (1,1)$, and $[f_{1},\dots,f_{9}]=[1,1,0,\dots,0]$, where $f_{1}=x_{11}x_{22}-x_{12}x_{21}$ and $f_{2}=x_{11}x_{33}-x_{13}x_{31}$. Similarly $\mathring{E_{2}} \cap \phi^{-1}(e) \cap \tilde{U}_{11}$ is given by
	\begin{flalign*}
		\eta=1,u_{12}=u_{13}=u_{21}=u_{31}=\theta=v_{23}'=v_{32}'=0,v_{22}'=v_{33}'=1.
	\end{flalign*}
	This yields $\chi(\mathring{E_{1}} \cap \phi^{-1}(e))=1$ and $Z^{\mathrm{mon}}_{Z_{2},e}=1-t$. 
	
	In conclusion, the possible values of $Z^{\mathrm{mon}}_{Z_{2},e}(t)$ are $1$, $1-t$ and $1-t^{2}$, verifying Theorem \ref{thmB}. 
	\end{example}
	
\section{Proof of the holomorphy conjecture}\label{sec5}
	In this section, we will introduce the holomorphy conjecture and prove it for determinantal varieties. Before that, we give the definition of twisted topological zeta function of an ideal associated to an integer.
	\begin{definition}
	Suppose $I \subset k[x_{1},\dots,x_{n}]$ is an ideal and $d \in \mathbb{Z}_{\ge 1}$ is an integer. Let $\mu:Y \rightarrow \mathbb{A}_{k}^{n}$ be a log resolution of $(\mathbb{A}_{k}^{n},V(I))$. Assume $E_{i}(i \in S)$ are the irreducible components of $\mu^{*}(I)$ and relative canonical divisor $K_{Y/\mathbb{A}_{k}^{n}}$ such that $\mu^{*}(I)=\sum_{i \in S}N_{i}E_{i}$ and $K_{Y/\mathbb{A}_{k}^{n}}=\sum_{i \in S}(\nu_{i}-1)E_{i}$. For any non-empty subset $J \subset S$, let $m_{J}:=\mathop{gcd}\limits_{i \in J}(N_{i})$. The twisted topological zeta function of $I$ associated to $d$ is defined to be
	\begin{flalign*}
		Z^{\mathrm{top},(d)}_{I}(s)=\sum_{J \subset S,d | m_{J}}\chi(\mathring{E}_{J})\prod_{i \in J}\frac{1}{N_{i}s+\nu_{i}},
	\end{flalign*}
	where $\mathring{E}_{J}:=(\cap_{i \in J}E_{i})\setminus (\cup_{j \notin J}E_{j})$. 
	\end{definition}
	
	\begin{remark}
	In \cite{The_holomorphy_conjecture_for_ideals_in_dimension_two}, the authors proved that the local version of this definition is independent of the choice of the resolution. The same proof shows that the global definition given above is also independent.
	\end{remark}
	
	The holomorphy conjecture is originally stated for a single polynomial as follows.
	\begin{conjecture}[Holomorphy conjecture, \cite{Geometry_on_arc_spaces_of_algebraic_varieties}]
	Let $f \in \mathbb{C}[x_{1},...,x_{n}]\setminus \mathbb{C}$, then $Z^{\mathrm{top},(d)}_{(f)}(s)$ is a polynomial, unless there exists an eigenvalue of the monodromy action such that its order is divisible by $d$.  
	\end{conjecture}
	
	This conjecture is an analogue of the monodromy conjecture and it can be generalized to the ideal case by replacing $(f)$ with an ideal $I$ and replacing Milnor monodromy with Verdier monodromy, see \cite{The_holomorphy_conjecture_for_ideals_in_dimension_two} for details. Similar to the monodromy conjecture, there are only some partial results about the conjecture. In \cite{Loeser_holo_conj_for_curve} and \cite{Veys_holomorphy_conj_for_curve}, Loeser and Veys proved the conjecture for two-dimensional principal ideal case. The general ideal case in dimension two was proved in \cite{The_holomorphy_conjecture_for_ideals_in_dimension_two}. In \cite{holomorphy_conj_for_nondegenerate} the authors solve the case of non-degenerate surface singularities.
	 
	From now on we use the setting in Scetion \ref{4} and we calculate the twisted topological zeta function of $I_{r}$ associated to any positive integer $d$ as in the following theorem.
	\begin{theorem}\label{twisted_top_zeta_function}
	The twisted topological zeta function of $I_{r}$ associated to $d$ is
	\begin{flalign*}
		Z^{\mathrm{top},(d)}_{I_{r}}(s)=\begin{cases}
			0,  \; d \ge 2, \\
			\prod_{j=1}^{r}\frac{1}{1+s \cdot \frac{r+1-j}{(m+1-j)(n+1-j)}}, \; d=1.
		\end{cases}
	\end{flalign*}
	\end{theorem}  
	\begin{proof}
		It suffices to compute $\chi(\mathring{E}_{I})$ for any $I \subset S=\{1, \dots, r\}$. By Proposition \ref{lemmabudur} and Proposition \ref{lemmadefloeser}, we can take $u_{i}$ to be $0$ or $1$ such that $I_{u_{1},\dots,u_{r}}=I$, and we have
	\begin{flalign*}
		[\mathring{E}_{I}]=[\mathfrak{X}^{l}_{u_{1},\dots,u_{r}}] \cdot \mathbb{L}^{-mnl+|I|+\sum_{i \in I}\nu_{i}-1} \cdot (\mathbb{L}-1)^{-|I|},
	\end{flalign*}
	where $l$ is an integer large enough such that Proposition \ref{lemmabudur} and Proposition \ref{lemmadefloeser} hold. Note that since the strict transform of $Z_{j}$ in $Y$ is $jE_{1}+(j-1)E_{2}+ \dots +E_{j}$ for all $1 \le j \le r$, $\mathfrak{X}^{l}_{u_{1},\dots,u_{r}}$ is the union of all $\mathcal{C}_{\lambda}^{l}$ such that $\lambda$ satisfies $\lambda_{j}=u_{1}+\dots+u_{j}$ for all $1 \le j \le r$. Suppose $\lambda_{r+1} \le \dots \le \lambda_{s} \le l$ and $\lambda_{s+1}=\dots=\lambda_{m}=\infty$, where $r \le s \le m$. Let \begin{flalign*}
		\lambda_{1}&=a_{1}, \\
		\lambda_{2}&=a_{1}+a_{2}, \\
		&\vdots \\
		\lambda_{s}&=a_{1}+\dots+a_{s},
	\end{flalign*}
	and we define $J:=\{2 \le i \le s \; | \; a_{i} \neq 0\}$. Suppose $J \cup \{1\}=\{i_{1},\dots,i_{p}\}$, where $i_{1}>\dots>i_{p}$, then by Proposition \ref{c_lambda} we have
	\begin{flalign*}
		[\mathcal{C}_{\lambda}^{l}]=\frac{[\mathrm{GL}_{m}] \cdot [\mathrm{GL}_{n}] \cdot \prod_{j=1}^{p}[G(i_{j-1}-i_{j},s+1-i_{j})]^{2} \cdot [\mathrm{GL}_{i_{j-1}-i_{j}}]}{[\mathrm{GL}_{s}]^{2} \cdot [\mathrm{GL}_{m-s}] \cdot [\mathrm{GL}_{n-s}] \cdot \mathbb{L}^{s^{2}l-(m+n)sl+s(m+n-2s)+\sum_{j=1}^{s}(m+n+1-2j)\lambda_{j}}},
	\end{flalign*}
	where we define $i_{0}=s+1$. Using the value given in Lemma \ref{Grassmannian}, the power of $(\mathbb{L}-1)$ in $[\mathcal{C}^{l}_{\lambda}]$ is $s$. Since $s \ge r \ge |I|$, $\chi(\mathring{E}_{I})=0$ unless $s=r=|I|$. In that case, we have $u_{1}=\dots=u_{r}=1$ and $\mathfrak{X}_{u_{1},\dots,u_{r}}=\mathcal{C}_{\lambda}^{l}$ where $\lambda_{j}=j$ for $1 \le j \le r$ and $\lambda_{j}=\infty$ for $j \ge r+1$, so $p=r$ and $i_{j}=r+1-j$ for $1 \le j \le r$. Now we obtain the value of $\chi(\mathring{E}_{I})$ as follows.
	\begin{align}
		\chi(\mathring{E}_{I})=&\frac{m! \cdot n! \cdot (r!)^{2}}{(r!)^{2} \cdot (m-r)! \cdot (n-r)!}  \\
		=&\frac{m! \cdot n!}{(m-r)! \cdot (n-r)!}.
	\end{align} 
	Since $N_{j}=r+1-j$ for all $1 \le j \le r$, for any integer $d \ge 2$, the topological zeta function of $Z_{r}$ associated to $d$ is $0$. The theorem now follows from direct computation using $(7)$ and the value of each $N_{j}, \nu_{j}$.
	\end{proof}
	
	\begin{remark}
	When $d=1$ and $m=n$, our result recovers Theorem \ref{roi_topological_zeta_function}. Although Docampo's work effectively utilizes orbit decomposition to compute the topological zeta function, this technique does not readily generalize to the twisted topological zeta function.
	\end{remark}

	Now we can prove the holomorphy conjecture for $I_{r}$.
	\begin{theorem}\label{thmC}
	The holomorphy conjecture holds for $I_{r}$.	
	\end{theorem}
	\begin{proof}
	By Theorem \ref{twisted_top_zeta_function}, the twisted topological zeta function of $I_{r}$ associated to an integer $d$ is a polynomial if and only if $d \ge 2$, so the conjecture holds directly.
	\end{proof}

	\end{document}